\documentclass[a4paper]{article}
\usepackage{fullpage}
\usepackage[centertags]{amsmath}
\usepackage{amsfonts}
\usepackage{amssymb}
\usepackage{amsthm}
\usepackage{amsmath}
\usepackage{cite} 

\newtheorem{theorem}{Theorem}

\newtheorem{lemma}[theorem]{Lemma}

\theoremstyle{remark}
\newtheorem{remark}[theorem]{Remark}
\theoremstyle{definition}
\newtheorem{definition}[theorem]{Definition}

\numberwithin{equation}{section}
\numberwithin{theorem}{section}

\newcommand{\Real}  {\mathbb R}
\newcommand{\V}     {\mathcal V}
\newcommand{\eps}{\varepsilon}

\newcommand{\norm}[1]{\left\Vert#1\right\Vert}
\newcommand{\abs}[1]{\left\vert#1\right\vert}
\newcommand{\lnorm}[1]{\vert#1\vert} 
\newcommand{\vnorm}[1]{\Vert#1\Vert} 

\title{Analytical Study of Certain Magnetohydrodynamic-$\alpha$ Models}
\author{Jasmine S.~Linshiz$^{1a}$
 and Edriss S.~Titi$^{1,2}$}
\date{} 

\begin{document}

\maketitle

\begin{center}
$^1$\textit{Department of Computer Science and Applied Mathematics \\
Weizmann Institute of Science  \\
Rehovot 76100, Israel}\\
$^a$jasmine.tal@weizmann.ac.il \\
$^2$\textit{Department of Mathematics
and\\ Department of Mechanical and  Aerospace Engineering \\
University of California \\
Irvine, CA  92697-3875, USA} \\
etiti@math.uci.edu \textit{and} edriss.titi@weizmann.ac.il
\end{center}

\begin{abstract}
In this paper we present an analytical study of a subgrid scale
turbulence model of the three-dimensional magnetohydrodynamic (MHD)
equations, inspired by the Navier-Stokes-$\alpha $ (also known as
the viscous Camassa-Holm equations or the Lagrangian-averaged
Navier-Stokes-$\alpha $ model. Specifically, we show the global
well-posedness and regularity of solutions of a certain MHD-$\alpha$
model (which is a particular case of the Lagrangian averaged
magnetohydrodynamic-$\alpha$  model without enhancing the
dissipation for the magnetic field).
We also introduce other subgrid scale turbulence models, inspired by
the Leray-$\alpha $ and the modified Leray-$\alpha $ models of
turbulence.
Finally, we discuss the relation of the  MHD-$\alpha $ model to the
MHD equations by proving a convergence theorem, that is, as the
length scale $\alpha $ tends to zero, a subsequence of solutions of
the MHD-$\alpha $ equations converges to a certain solution (a
Leray-Hopf solution) of the three-dimensional MHD equations.
\end{abstract}

\textbf{Keywords:} subgrid scale models; turbulence models;
magnetohydrodynamics; regularizing MHD; magnetohydrodynamic-$ \alpha
$ model; Lagrangian-averaged magnetohydrodynamic-$ \alpha $ model;
Leray-$\alpha$ model.

\textbf{Mathematics Subject Classification:}
%
76D03, 
76F20, 
76F55, 
76F65,  
76W05.  

\section{Introduction}
We consider the three-dimensional magnetohydrodynamic (MHD)
equations for a homogeneous incompressible resistive viscous fluid
subjected to a Lorentz force due to the presence of a magnetic
field. The MHD involves coupling Maxwell's equations governing the
magnetic field and the Navier-Stokes equations (NSE) governing the
fluid motion. The system has the form
\begin{align}\label{grp:MHD}
\begin{split}
& \frac{\partial \boldsymbol{v}}{\partial t}+\left(
\boldsymbol{v}\cdot \nabla \right) \boldsymbol{v}-\nu \Delta
\boldsymbol{v}+\nabla \pi + \frac{1}{2}\nabla \lnorm{
\boldsymbol{B}} ^{2}=\left( \boldsymbol{B}\cdot \nabla \right)
\boldsymbol{B},
\\
& \frac{\partial \boldsymbol{B}}{\partial t}+\left(
\boldsymbol{v}\cdot \nabla \right) \boldsymbol{B}-\left(
\boldsymbol{B}\cdot \nabla \right) \boldsymbol{v}-\eta \Delta
\boldsymbol{B} =0,
\\
& \nabla \cdot \boldsymbol{v}=\nabla \cdot \boldsymbol{B}=0,
\end{split}
\end{align}
where $\boldsymbol{v}\left( x,t\right) $, the fluid velocity field,
$ \boldsymbol{B}\left( x,t\right) $, the magnetic field and $\pi$,
the pressure, are the unknowns; $\nu >0$ is the constant kinematic
viscosity and $\eta >0$ is the constant magnetic diffusivity.

Current scientific methods and tools are unable to compute the
turbulent behavior of three-dimensional (3D) fluids and
magnetofluids analytically or via direct numerical simulation due to
the large range of scales of motion that need to be resolved when
the Reynolds number is high. For many purposes, it might be adequate
to compute only certain statistical features of the physical
phenomenon of turbulence and much effort is being made to produce
reliable turbulence models that parameterize  the average effects of
the
fluctuations on the averages, without calculating the former explicitly. %
Motivated by the remarkable performance of the
Navier-Stokes-$\alpha$ (NS-$ \alpha $) (also known as the viscous
Camassa-Holm equations (VCHE) or the Lagrangian-averaged
Navier-Stokes-$\alpha $ (LANS-$\alpha $)) as a closure model of
turbulence in infinite channels and pipes, whose solutions give
excellent agreement with empirical data for a wide range of large
Reynolds numbers \cite{a_CFHOTW99, a_CFHOTW99_ChanPipe,a_CFHOTW98},
the alpha subgrid scale models of turbulence have been extensively
studied in recent years (see, e.g.,
\cite{a_HT05,a_CHOT05,a_ILT05,a_VTC05,a_CTV05,a_FHT02, a_FHT01,
a_CFHOTW99, a_CFHOTW99_ChanPipe,a_CFHOTW98,
a_MKSM03,a_LL06,a_LL03,a_L06}).
%

A justification of the inviscid NS-$\alpha $ model can be found, for
example, in
\cite{a_CFHOTW99_ChanPipe,a_HMR98,a_H02_pA,a_MS03,a_C01}.

An extension of the NS-$ \alpha $ model to the nondissipative MHD
is given, e.g., in \cite{a_H02_ch}. 
 The model was obtained
from variational principles by modifying the Hamiltonian associated
with the ideal MHD equations subject to the incompressibility
constraint. Then the dissipation is introduced in an \textit{ad hoc}
fashion in analogy to the NS-$\alpha $, following
\cite{a_CFHOTW98,a_CFHOTW99,a_CFHOTW99_ChanPipe,a_FHT02}.
Specifically, the flow Lagrangian of the ideal MHD is given by
\begin{equation*}
\mathcal L[\boldsymbol u,D,\boldsymbol B]= \int \left(\frac{1}{2} D
|\boldsymbol u|^2 - \pi
 (D-1) - \frac{1}{2} |\boldsymbol
B|^2\right)dx
\end{equation*}
with volume preservation for the pressure.  Here the volume element
$D(x,t)=\left(\det \left({\partial {X}}/{\partial
a}\right)(a,t)\right)^{-1}$ at \mbox{$x=X(a,t)$}, where
\mbox{${X}(a,t)$} is the Lagrangian fluid trajectory, ${\partial
{X}}/{\partial t}(a,t)=\boldsymbol{u}(x,t)$ (see \cite{a_H99}).
First, the Lagrangian is averaged and approximated using a form of
Taylor's hypothesis (see, e.g., \cite{a_H02_fD}) to obtain
\begin{equation*}
\bar{\mathcal L }= \int \left(\frac{1}{2} D \left(|\boldsymbol
u|^2+\alpha^2|\nabla \boldsymbol u|^2\right) - \pi (D-1) -
\frac{1}{2} \left(|\boldsymbol B|^2+\alpha_M^2|\nabla \boldsymbol
B|^2\right)\right)dx,
\end{equation*}
then the Hamiltonian principle is applied (see, e.g.,
\cite{a_HMR98}) to produce an ideal MHD-$\alpha $ model (eq.
\eqref{grp:LAMHD} with $\nu=\eta=0$). Adding viscosity  and
diffusivity provides the MHD-$\alpha $ (or the Lagrangian-averaged
magnetohydrodynamic-$ \alpha $ (LAMHD-$ \alpha $)) model
\begin{align}\label{grp:LAMHD}
\begin{split}
& \frac{\partial \boldsymbol{v}}{\partial t}+\left(
\boldsymbol{u}\cdot \nabla \right)
\boldsymbol{v}+\sum_{j=1}^{3}\boldsymbol{v}_{j}\nabla
\boldsymbol{u}_{j}-\nu \Delta \boldsymbol{v}+\nabla
p+%
\sum_{j=1}^{3}(\boldsymbol{B_s})_{j}\nabla \boldsymbol{B}_{j}=\left(
\boldsymbol{B_s}\cdot \nabla \right) \boldsymbol{B},
\\
& \frac{\partial \boldsymbol{B_s}}{\partial t}+\left(
\boldsymbol{u}\cdot \nabla \right) \boldsymbol{B_s}-\left(
\boldsymbol{B_s}\cdot \nabla \right) \boldsymbol{u}-\eta \Delta
\boldsymbol{B}%
=0 ,
\\
& \boldsymbol{v}=\left( 1-\alpha ^{2}\Delta \right) \boldsymbol{u},
\qquad  \boldsymbol{B}=\left( 1-\alpha_M
^{2}\Delta \right) \boldsymbol{B_s}, \\
& \nabla \cdot \boldsymbol{u}=\nabla \cdot \boldsymbol{v}=\nabla
\cdot \boldsymbol{B_s}=\nabla \cdot \boldsymbol{B}=0,
\end{split}
\end{align}
where $\boldsymbol{u}$ and $\boldsymbol{B_s}$ represent the unknown
`filtered' fluid velocity  and magnetic fields, respectively, $p$ is
the unknown `filtered' pressure,  and $\alpha
>0, \, \alpha_M>0$ are lengthscale parameters that represent the width of the
filters. At the limit $\alpha =0,\,\alpha_M=0$, we formally obtain
the three-dimensional MHD equations. The LAMHD-$ \alpha $ model was
investigated numerically in periodic boundary conditions in two
\cite{a_MMP05_2D,a_GHMP06} and three \cite{a_MMP05_3D} space
dimensions against direct numerical simulations. %
In \cite{a_GHMP06} the K\'arm\'an-Howarth theorem was extended to
LAMHD-$ \alpha $ equations. The LAMHD-$ \alpha $ model was also
studied in \cite{a_JR05} in the context of convection-driven plane
layer geodynamo models.

We tend to think about the $\alpha$ models as a numerical
regularization of the underlying equation, which makes the
nonlinearity milder, and hence the solutions of the modified
equation are smoother. This is contrary to the hyperviscosity
regularization \cite{a_L59}  and nonlinear viscosity
\cite{a_L70,b_L85,a_S63}, which lead to unnecessary extra
dissipation of the energy of the system. To emphasize this numerical
analysis point of view, we observe that recently a Leray-$\alpha$
model of the inviscid Burgers equation
\begin{equation}\label{eq:Burgers}
\frac{\partial \boldsymbol{v}}{\partial t} + \boldsymbol{v}
\frac{\partial \boldsymbol{v}}{\partial x} =0,
\end{equation}
which is
\begin{align}\label{eq:BurgersAlpha}
\begin{cases}
&\frac{\partial \boldsymbol{v}^\alpha}{\partial t} +
\boldsymbol{u}^\alpha
\frac{\partial \boldsymbol{v}^\alpha}{\partial x} =0, \\
& \boldsymbol{v}^\alpha=\boldsymbol{u}^\alpha-\alpha
^{2}\boldsymbol{u}^\alpha_{xx},
\end{cases}
\end{align}
has been introduced in \cite{A_BF06} and \cite{a_TTZ06}. Regular
unique solutions of \eqref{eq:BurgersAlpha} exist globally and it
was shown computationally in \cite{A_BF06} that the solutions of
\eqref{eq:BurgersAlpha} converge to the unique entropy weak solution
(see, e.g., \cite{a_O63,b_S94,a_TTZ06}) of \eqref{eq:Burgers}.
Notice that there is no dissipation in \eqref{eq:BurgersAlpha}, and
the $L^\infty$ norm of $\boldsymbol{v}^\alpha$ is preserved. On the
other hand, the viscous regularizing approach, which is usually
taken for the Burgers equation, is achieved by introducing an
artificial viscosity term in \eqref{eq:Burgers} and obtaining the
viscous Burgers equation
\begin{equation}\label{eq:BurgersViscous}
\frac{\partial \boldsymbol{v}^\eps}{\partial t}-\eps^2
\frac{\partial^2 \boldsymbol{v}^\eps}{\partial x^2}+
\boldsymbol{v}^\eps \frac{\partial \boldsymbol{v}^\eps}{\partial x}
=0.
\end{equation}
This model gives a smooth solution $\boldsymbol{v}^{\eps}$, which
converges in the appropriate norms to the unique entropy weak
solution (see, e.g., \cite{a_O63}). However, the energy of
$\boldsymbol{v}^{\eps}$  is decaying in time at a much higher rate
than the decay expected for the entropy weak solution. Hence, the
advantage of introducing the Leray-$\alpha$ model
\eqref{eq:BurgersAlpha} for Burgers equation. This simple example
clarifies our numerical approach of why we insist on making the
nonlinearity milder instead of adding additional viscous or
hyperviscous terms. This approach has been discussed further in
\cite{a_CLT06} in the context of Euler and Navier-Stokes equations.

Filtering the magnetic field, as it is done in
\cite{a_MMP05_3D,a_MMP05_2D,a_JR05}, is equivalent to introducing
hyperdiffusivity for the filtered magnetic field $\boldsymbol{B_s}$,
due to the term $-\eta\alpha_M^2 \Delta^2 \boldsymbol{B_s}$ in
\eqref{grp:LAMHD}, which we think is unnecessary. Taking the
numerical analysis point of view discussed above we prove the
well-posedness of a certain MHD-$\alpha$ model without introducing
extra dissipation for the magnetic field, i.e.~we filter only the
velocity field, but not the magnetic field and obtain the following
regularizing system of \eqref{grp:MHD}
\begin{align}\label{grp:alphaMHD_intro}
\begin{split}
& \frac{\partial \boldsymbol{v}}{\partial t}+\left(
\boldsymbol{u}\cdot \nabla \right)
\boldsymbol{v}+\sum_{j=1}^{3}\boldsymbol{v}_{j}\nabla
\boldsymbol{u}_{j}-\nu \Delta \boldsymbol{v}+\nabla
p+\frac{1}{2}\nabla \lnorm{ \boldsymbol{B}} ^{2}=\left(
\boldsymbol{B}\cdot \nabla
\right) \boldsymbol{B},   \\
& \frac{\partial \boldsymbol{B}}{\partial t}+\left(
\boldsymbol{u}\cdot \nabla \right) \boldsymbol{B}-\left(
\boldsymbol{B}\cdot \nabla \right) \boldsymbol{u}-\eta \Delta
\boldsymbol{B}%
=0 ,\\
& \boldsymbol{v}=\left( 1-\alpha ^{2}\Delta \right) \boldsymbol{u}, \qquad \alpha>0,\\
& \nabla \cdot \boldsymbol{u}=\nabla \cdot \boldsymbol{v}=\nabla
\cdot \boldsymbol{B}=0,
\end{split}
\end{align}
instead of the system \eqref{grp:LAMHD}.

As $\alpha$ models are some sort of regularizing numerical schemes,
we would like to make sure that they inherit some of the original
properties of the 3D MHD equations.
Formally, three ideal, i.e.~$\nu=\eta  =0$, quadratic invariants of
the system \eqref{grp:alphaMHD_intro} could be identified with the
invariants of the original ideal 3D MHD equations under suitable
boundary conditions, for instance, in rectangular periodic boundary
conditions or in the whole space $\Real^3$. Namely, the energy
\mbox{$E^{\alpha }=\frac{1}{2}\int_{\Omega }\left(
\boldsymbol{v}\left( x\right) \cdot \boldsymbol{u}(x)+\lnorm{
\boldsymbol{B}(x)} ^{2}\right) dx$}, the cross helicity \mbox{$
H_{C}^{\alpha }=\frac{1}{2}\int_{\Omega }\boldsymbol{v}(x)\cdot
\boldsymbol{B }(x)dx$}, and the magnetic helicity \mbox{$
H_{M}^{\alpha }=\frac{1}{2}\int_{\Omega }\boldsymbol{A}(x)\cdot
\boldsymbol{B }(x)dx$}, where $ \boldsymbol{A}$ is the vector
potential, so that \mbox{$\boldsymbol{B}= \nabla \times
\boldsymbol{A}$}; and they reduce, as $\alpha \rightarrow 0$, to the
ideal invariants of the MHD equations.

There are other possible alpha subgrid scale models that can be
shown to have global existence and uniqueness. For instance,
inspired by the Leray-$\alpha $
\cite{a_CHOT05,a_CTV05,a_VTC05,a_HN03,a_GH03} and modified
Leray-$\alpha $ \mbox{(ML-$ \alpha $)} \cite{a_ILT05} models of
turbulence, we formulate similar MHD alpha models, we refer to them
as Leray-$\alpha $-MHD and ML-$ \alpha $-MHD models, respectively.
The Leray-$\alpha $ and ML-$ \alpha $ models of turbulence reduce to
the same closure model for the Reynolds averaged Navier-Stokes
equations in turbulent channels and pipes as the NS-$\alpha $ model
under the corresponding symmetries \cite{a_CHOT05,a_CTV05,a_ILT05},
which, as we mentioned above, compares successfully with
experimental data for a wide range of Reynolds numbers. This
comparison means that the Leray-$\alpha $ and the ML-$\alpha $
models as well as NS-$\alpha $ equations could be equally used as
subgrid scale models of turbulence.

Specifically, we consider the following version of the three-dimensional Leray-$\alpha $%
-MHD model 
\begin{align}\label{grp:Leray_alpha_MHD_intro}
\begin{split}
& \frac{\partial \boldsymbol{v}}{\partial t}+\left(
\boldsymbol{u}\cdot \nabla \right) \boldsymbol{v}-\nu \Delta
\boldsymbol{v}+\nabla p+\frac{1}{ 2}\nabla \lnorm{ \boldsymbol{B}}
^{2}=\left( \boldsymbol{B}
\cdot \nabla \right) \boldsymbol{B},   \\
& \frac{\partial \boldsymbol{B}}{\partial t}+\left(
\boldsymbol{u}\cdot \nabla \right) \boldsymbol{B}-\left(
\boldsymbol{B}\cdot \nabla \right) \boldsymbol{v}-\eta \Delta
\boldsymbol{B}
 =0,
\\
& \boldsymbol{v}=\left( 1-\alpha ^{2}\Delta \right) \boldsymbol{u}, \\
& \nabla \cdot \boldsymbol{u}=\nabla \cdot \boldsymbol{v}=\nabla
\cdot \boldsymbol{B}=0.
\end{split}
\end{align}
Formally, the term $\left( \boldsymbol{B}\cdot \nabla \right)
\boldsymbol{v }$ comes from requiring in the ideal ($\nu=0=\eta$)
case the conservation of energy \mbox{$ E^{\alpha
}=\frac{1}{2}\int_{\Omega }\left( \lnorm{ \boldsymbol{v} (x)}
^{2}+\lnorm{ \boldsymbol{B}(x)} ^{2}\right) dx$} (under suitable
boundary conditions). While the requirement for the system to have
an ideal invariant corresponding to the cross helicity \mbox{$
H_{C}^{\alpha }=\frac{1}{2}\int_{\Omega }\boldsymbol{v}(x)\cdot
\boldsymbol{B }(x)dx  $} leads to the term $\left(
\boldsymbol{u}\cdot \nabla \right) \boldsymbol{B}$.
Contrary to the MHD-$\alpha $ model \eqref{grp:alphaMHD_intro},
where we establish the
 existence and uniqueness, for the 3D Leray-$\alpha $-MHD model
 \eqref{grp:Leray_alpha_MHD_intro} we are able to
establish only the existence of weak solutions, as in the case for
the original MHD equations \eqref{grp:MHD}.
In this case, the term $\left( \boldsymbol{B}\cdot \nabla \right)
\boldsymbol{v}$ is problematic as in the usual 3D NSE and MHD.
However, in the two dimensional case the existence and uniqueness of
weak solutions can be shown (similarly to the proof given for the
model \eqref{grp:alphaMHD_intro} in section 3) for the following
2D-Leray-$\alpha $-MHD model
\begin{align}\label{grp:2D_Leray_alpha_MHD}
\begin{split}
& \frac{\partial \boldsymbol{v}}{\partial t}+\left(
\boldsymbol{u}\cdot
\nabla \right) \boldsymbol{v}-\nu \Delta \boldsymbol{v}+\nabla p+\frac{1}{%
2}\nabla \lnorm{ \boldsymbol{B}} ^{2}=\left( \boldsymbol{B}
\cdot \nabla \right) \boldsymbol{B}, \\
& \frac{\partial \boldsymbol{B}}{\partial t}+\left(
\boldsymbol{u}\cdot \nabla \right) \boldsymbol{B}-\left(
\boldsymbol{B}\cdot \nabla \right)
\boldsymbol{u}-\eta \Delta \boldsymbol{B}%
=0, \\
& \boldsymbol{v}=\left( 1-\alpha ^{2}\Delta \right) \boldsymbol{u}, \\
& \nabla \cdot \boldsymbol{u}=\nabla \cdot \boldsymbol{v}=\nabla
\cdot \boldsymbol{B}=0.
\end{split}
\end{align}
For this system, due to the identity \mbox{$\int_{\Omega }\left(
\boldsymbol{u}\cdot \nabla \boldsymbol{u}\right) \cdot \Delta
\boldsymbol{u}=0
$} %
(for the periodic 2D case and divergence free $\boldsymbol{u}$), the
ideal invariant corresponding to the energy is \mbox{$ E^{\alpha
}=\frac{1}{2}\int_{\Omega }\left( \boldsymbol{v}\left( x\right)
\boldsymbol{u}(x)+\lnorm{ \boldsymbol{B}(x)} ^{2}\right) dx$}. %
At the moment we are unable to find a conserved quantity in the
ideal version of \eqref{grp:2D_Leray_alpha_MHD} that can be
identified with a cross helicity. The mean-square magnetic
potential, given by $\mathcal{A}=\frac{1}{2}\int_{\Omega } \lnorm{
\psi(x)} ^{2} dx$, where $\boldsymbol{B}=\nabla^{\perp}\psi$, is
conserved in the ideal case. We note that it appears that there is
no conserve quantity that could be identified with energy for the 3D
version of \eqref{grp:2D_Leray_alpha_MHD}.

The three-dimensional Modified-Leray-$ \alpha $-MHD model, for which
the well-posedness can be proved in a similar way as for the model
\eqref{grp:alphaMHD_intro},
is given by
\begin{align}
\begin{split}\label{grp:ML_alpha_MHD_intro}
& \frac{\partial \boldsymbol{v}}{\partial t}+\left(
\boldsymbol{v}\cdot
\nabla \right) \boldsymbol{u}-\nu \Delta \boldsymbol{v}+\nabla p+\frac{1}{%
2}\nabla \lnorm{ \boldsymbol{B}} ^{2}=\left( \boldsymbol{B}%
\cdot \nabla \right) \boldsymbol{B},  \\
& \frac{\partial \boldsymbol{B}}{\partial t}+\left(
\boldsymbol{u}\cdot \nabla \right) \boldsymbol{B}-\left(
\boldsymbol{B}\cdot \nabla \right) \boldsymbol{u}-\eta \Delta
\boldsymbol{B}%
=0, \\
& \boldsymbol{v}=\left( 1-\alpha ^{2}\Delta \right) \boldsymbol{u}, \\
& \nabla \cdot \boldsymbol{u}=\nabla \cdot \boldsymbol{v}=\nabla
\cdot \boldsymbol{B}=0,
\end{split}
\end{align}
where the term $\left( \boldsymbol{B}\cdot \nabla \right)
\boldsymbol{u}$ comes from requiring the conservation of energy (in
the ideal case, with periodic boundary conditions or in $\Real^3$)
\mbox{$ E^{\alpha }=\frac{1}{2}\int_{\Omega }\left(
\boldsymbol{v}\left( x\right) \boldsymbol{u}(x)+\lnorm{
\boldsymbol{B}(x)} ^{2}\right) dx$}. Also, the system conserves the
magnetic helicity \mbox{$ H_{M}^{\alpha }=\frac{1}{2}\int_{\Omega
}\boldsymbol{A}(x)\cdot \boldsymbol{B }(x)dx$}. At the moment we are
unable to find a conserved quantity in the ideal version of
\eqref{grp:ML_alpha_MHD_intro} which can be identified with a cross
helicity.

The main goal of this paper is to establish the global existence,
uniqueness and regularity of solutions of the three-dimensional
MHD-$\alpha $ equations \eqref{grp:alphaMHD_intro} subject to
periodic boundary conditions (similar results also hold in $\Real
^3$). We emphasize again that we consider a version of the MHD alpha
models, where only the velocity field is filtered, while the
magnetic field remains unfiltered. We note that in the case of
filtering the magnetic field, as in \eqref{grp:LAMHD}, one has
hyperdiffusivity for the filtered  magnetic field $\boldsymbol{B_s}$
and the proof of the existence and uniqueness of regular solutions
of \eqref{grp:LAMHD} is deduced in a similar way.

We start by introducing some preliminary background and the
functional setting in section 2. %
In section 3 we show the global well-posedness of the MHD-$\alpha $
subgrid scale model of turbulence \eqref{grp:alphaMHD_intro}. We
remark that using the Gevrey regularity techniques developed in
\cite{a_FT89} (see also \cite{a_FT98}) one can show that the
solution of the MHD-$\alpha$ model becomes instantaneously analytic
in space and time. As a result of this Gevrey regularity, one
deduces the existence of a dissipation range in the energy spectrum
in which the energy decays exponentially fast as a function of the
wavenumber, for $k$ larger than the dissipation length scale (see
\cite{a_DT95}). One can also establish, in the forced case, the
existence of a finite dimensional global attractor, a subject of
future work.
In section 4 we relate the solutions of the MHD-$\alpha $ equations
to those of the 3D MHD as the length scale $\alpha $ tends to zero.
Specifically, we prove that one can extract subsequences of weak
solutions of the MHD-$ \alpha $ equations which converge as $\alpha
$ $\rightarrow $ $0^{+}$ (in the appropriate sense) to a Leray-Hopf
weak solution of the three-dimensional  MHD equations
\eqref{grp:MHD} on any time interval $[0,T]$, which satisfies the
energy inequality
\begin{equation*}
\lnorm{ \boldsymbol{v}\left( t\right) } ^{2}+\lnorm{ \boldsymbol{B}
\left( t\right) } ^{2} +2 \int_{t_0}^{t}\left( \nu\vnorm{
\boldsymbol{v}(s)} ^{2}+\eta\vnorm{ \boldsymbol{B} (s)}
^{2}\right)ds \leq \lnorm{ \boldsymbol{v}\left( t_0\right) }
^{2}+\lnorm{ \boldsymbol{B}\left( t_0\right) } ^{2}
\end{equation*}
for almost every $t_0\in[0,T]$ and all $t\in\left[t_0,T\right]$.
Also, if the initial data is smooth a subsequence of solutions
converges for a short interval of time, that depends on the initial
data, $\nu$, $\eta$ and the domain, to the unique strong solution of
the MHD equations on this interval. Thus the $\alpha$ models can be
viewed as a regularizing numerical method. Section 5 contains a
discussion summarizing our results.


\section{Functional Setting and Preliminaries}

Let $\Omega $ be the $L$-periodic three-dimensional box $\Omega
=[0,L]^{3}$. We consider the following  MHD-$\alpha $ subgrid scale
turbulence model, which we introduced in \eqref{grp:alphaMHD_intro},
subject to periodic boundary condition with a basic domain $\Omega
$,
\begin{subequations}
\label{grp:alphaMHD}
\begin{align}
& \frac{\partial \boldsymbol{v}}{\partial t}+\left(
\boldsymbol{u}\cdot \nabla \right)
\boldsymbol{v}+\sum_{j=1}^{3}\boldsymbol{v}_{j}\nabla
\boldsymbol{u}_{j}-\nu \Delta \boldsymbol{v}+\nabla
p+\frac{1}{2}\nabla \lnorm{ \boldsymbol{B}} ^{2}=\left(
\boldsymbol{B}\cdot \nabla
\right) \boldsymbol{B},  \label{eq:alphaMHD:velocity} \\
& \frac{\partial \boldsymbol{B}}{\partial t}+\left(
\boldsymbol{u}\cdot \nabla \right) \boldsymbol{B}-\left(
\boldsymbol{B}\cdot \nabla \right) \boldsymbol{u}-\eta \Delta
\boldsymbol{B}%
=0,
\label{eq:alphaMHD:magneticField} \\
& \boldsymbol{v}=\left( 1-\alpha ^{2}\Delta \right) \boldsymbol{u}, \\
& \nabla \cdot \boldsymbol{u}=\nabla \cdot \boldsymbol{v}=\nabla
\cdot
\boldsymbol{B}=0,  \label{eq:alphaMHD:divFreeCond} \\
& \boldsymbol{u}(x,0)=\boldsymbol{u}^{in}(x), \\
& \boldsymbol{B}(x,0)=\boldsymbol{B}^{in}(x),
\end{align}
\end{subequations}
where $\boldsymbol{u}$ represents the unknown `filtered' fluid
velocity vector, $p$ is the unknown `filtered' pressure, and $\alpha
>0$ is a lengthscale parameter which represents the width of the
filter. At the limit $\alpha =0$ we formally obtain the
three-dimensional MHD equations \eqref{grp:MHD}, where
$\boldsymbol{u}$ is the Eulerian velocity field and
$p-\frac{1}{2}|\boldsymbol{u}|^{2}$ is the pressure. Notice that we
chose to smooth only the velocity field and not the magnetic field,
thus we do not introduce hyperdiffusivity for the magnetic field, as
it is for the filtered magnetic field in \eqref{grp:LAMHD}.

We consider initial values with zero spatial means, i.e., we assume
that
\begin{equation}
\int_{\Omega }\boldsymbol{u}^{in}dx=\int_{\Omega
}\boldsymbol{B}^{in}dx=0, \label{eq:zeroSpatialMean}
\end{equation}
then from \eqref{eq:alphaMHD:velocity} and
\eqref{eq:alphaMHD:magneticField}, after integration by parts, using
the spatial periodicity of the solution and the divergence free
condition \eqref{eq:alphaMHD:divFreeCond} we have
\mbox{$\left({d}/{dt}\right)\int_{\Omega }\boldsymbol{v}dx=0$},
\mbox{$\left({d}/{dt}\right) \int_{\Omega }\boldsymbol{B}dx=0$} and
\mbox{$\left({d}/{dt}\right)\int_{\Omega }\boldsymbol{ u}dx=0$}.
Namely, the spatial mean of the solution is invariant under time.
Hence, by \eqref{eq:zeroSpatialMean}, $\int_{\Omega }\boldsymbol{v}
dx=\int_{\Omega }\boldsymbol{u}dx=\int_{\Omega }\boldsymbol{B}dx=0$.

Next, we introduce some notation and background following the
mathematical theory of NSEs, see, for instance,
\cite{b_CF88,b_T95,b_T84,b_L85}. Let $ L^{p}(\Omega )$ and
$H^{m}(\Omega )$ denote the $L^{p}$ Lebesgue spaces and Sobolev
spaces respectively. We denote by $\lnorm{ \cdot } $ the
$L^{2}$-norm, and by $\left( \cdot ,\cdot \right) $ the
$L^{2}$-inner product. Let $X$ be a linear subspace of integrable
functions defined on the domain $\Omega $, we define
$\dot{X}:=\{\varphi \in X:\int_{\Omega }\varphi (x)dx=0\}$ and
{$\mathcal{V}=\{\varphi :\varphi \text{ is a vector valued
trigonometric polynomial defined on }\Omega ,\text{ such that
}\nabla \cdot \varphi =0\text{ and }\int_{\Omega }\varphi
(x)dx=0\}$}.
The spaces $H$ and $V$ are the closures of $\mathcal{V}$ in
$L^{2}(\Omega )$ and in $H^{1}(\Omega )$ respectively; observe that
$H^{\perp }$, the orthogonal complement of $H$ in $L^{2}(\Omega )$
is \{${\nabla p:p\in H^{1}(\Omega )}$\}. Let $P_{\sigma
}:\dot{L}^{2}\left( \Omega \right)
\rightarrow H$ be the Helmholtz-Leray projection, and \mbox{$%
A=-P_{\sigma }\Delta $} be the Stokes operator with domain $%
D(A)=(H^{2}(\Omega )\cap V)$. In the periodic boundary conditions
$A=-\Delta |_{D(A)}$ is a self-adjoint positive
operator with compact inverse. Hence the space $H$ has an orthonormal basis $%
\{w_{j}\}_{j=1}^{\infty }$ of eigenfunctions of $A$, $Aw_{j}=\lambda
_{j}w_{j}$, with $0<\lambda _{1}\leq \lambda _{2}\leq \ldots $,
$\lambda _{j}\sim j^{2/d}L^{-2}$, see, e.g., \cite{b_CF88,a_M78}.
One can show that $V=D\left( A^{1/2}\right) $. We denote
\mbox{$\left( \left( \cdot ,\cdot \right) \right) =\left(
A^{1/2}\cdot ,A^{1/2}\cdot \right) $} and $\vnorm{ \cdot } =\lnorm{
A^{1/2}\cdot } $ the inner product and the norm on $V$,
respectively.

Following the notation of the Navier-Stokes equations and those of
\cite {a_FHT02}, we denote
\begin{align*}
B\left( \boldsymbol{u},\boldsymbol{v}\right) & =P_{\sigma }\left[
\left(
\boldsymbol{u}\cdot \nabla \right) \boldsymbol{v}\right] ,\emph{\quad }%
\boldsymbol{u},\boldsymbol{v}\in \V, \\
\tilde{B}\left( \boldsymbol{u},\boldsymbol{v}\right) & =P_{\sigma
}\left[
\left( \nabla \times \boldsymbol{v}\right) \times \boldsymbol{u}\right] ,%
\emph{\quad }\boldsymbol{u},\boldsymbol{v}\in \V.
\end{align*}%
Notice that
\begin{equation*}
\left( B\left( \boldsymbol{u},\boldsymbol{v}\right)
,\boldsymbol{w}\right) =-\left( B\left(
\boldsymbol{u},\boldsymbol{w}\right) ,\boldsymbol{v}\right)
,\emph{\quad }\boldsymbol{u},\boldsymbol{v},\boldsymbol{w}\in \V,
\end{equation*}%
and due to the identity
\begin{equation}\label{eq:gen_3D_vector_id}
\left( b\cdot \nabla \right) a+\sum_{j=1}^{3}a_{j}\nabla
b_{j}=-b\times \left( \nabla \times a\right) +\nabla \left( a\cdot
b\right) ,
\end{equation}%
\begin{equation*}
\left( \tilde{B}\left( \boldsymbol{u},\boldsymbol{v}\right) ,\boldsymbol{w}%
\right) =\left( B\left( \boldsymbol{u},\boldsymbol{v}\right) ,\boldsymbol{w}%
\right) -\left( B\left( \boldsymbol{w},\boldsymbol{v}\right) ,\boldsymbol{u}%
\right) .
\end{equation*}%
The definitions of $B\left( \boldsymbol{u},\boldsymbol{v}\right) $
and $ \tilde{B}\left( \boldsymbol{u},\boldsymbol{v}\right) $ and the
above algebraic identities may be extended to larger spaces by the
density of $\V $ in the appropriate space each time the
corresponding trilinear forms are continuous. The extensions of the
bilinear forms $B$ and $\tilde{B}$ (which we also denote $B$ and
$\tilde{B}$) have the following properties

\begin{lemma}
\label{lemma:B_estimates} \mbox{}

\begin{enumerate}
\item \label{enu:lemma:BandBtilde_estimate_V_V_V}Let $X$ be either $B$ or $%
\tilde{B}$. The operator $X$ can be extended continuously from
$V\times V$ with values in $V^{\prime }$ (the dual space of $V$). In
particular, for every
$\boldsymbol{u},\boldsymbol{v},\boldsymbol{w}\in V$,
\begin{equation}
\abs{ \left\langle X\left( \boldsymbol{u},\boldsymbol{\boldsymbol{v}}%
\right) ,\boldsymbol{w}\right\rangle _{V^{\prime }}} \leq
c\lnorm{ \boldsymbol{u}} ^{1/2}\vnorm{ \boldsymbol{u}%
} ^{1/2}\vnorm{ \boldsymbol{\boldsymbol{v}}} \vnorm{ \boldsymbol{w}}
.  \label{eq:BandBtilde_estimate_V_V_V}
\end{equation}%
Moreover,
\begin{equation}
\left( B\left( \boldsymbol{u},\boldsymbol{v}\right)
,\boldsymbol{w}\right) =-\left( B\left(
\boldsymbol{u},\boldsymbol{w}\right) ,\boldsymbol{v}\right)
,\emph{\quad }\boldsymbol{u},\boldsymbol{v},\boldsymbol{w}\in V,
\label{eq:B_id1}
\end{equation}%
which in turn implies that
\begin{equation}
\left( B\left( \boldsymbol{u},\boldsymbol{v}\right)
,\boldsymbol{v}\right) =0,\emph{\quad
}\boldsymbol{u},\boldsymbol{v}\in V.  \label{eq:B_id2}
\end{equation}%
Also
\begin{equation}
\left( \tilde{B}\left( \boldsymbol{u},\boldsymbol{v}\right) ,\boldsymbol{w}%
\right) =\left( B\left( \boldsymbol{u},\boldsymbol{v}\right) ,\boldsymbol{w}%
\right) -\left( B\left( \boldsymbol{w},\boldsymbol{v}\right) ,\boldsymbol{u}%
\right) ,\emph{\quad
}\boldsymbol{u},\boldsymbol{v},\boldsymbol{w}\in V,
\label{eq:Btilda_id1}
\end{equation}%
and hence
\begin{equation}
\left( \tilde{B}\left( \boldsymbol{u},\boldsymbol{v}\right) ,\boldsymbol{u}%
\right) =0,\emph{\quad }\boldsymbol{u},\boldsymbol{v}\in V.
\label{eq:Btilda_id2}
\end{equation}

\item \label{enu:lemma:BandBtilde_estimate_DA_V_H}  Furthermore, let $\boldsymbol{u}%
\in D(A),\boldsymbol{v}\in V,\boldsymbol{w}\in H$ and let $X$ be
either $B$ or $\tilde{B}$ then
\begin{equation}
\abs{ \left( X\left( \boldsymbol{u},\boldsymbol{v}\right) ,%
\boldsymbol{w}\right) } \leq c\vnorm{ \boldsymbol{u}%
} ^{1/2}\lnorm{ A\boldsymbol{u}} ^{1/2}\vnorm{ \boldsymbol{v}}
\lnorm{ \boldsymbol{w}} . \label{eq:BandBtilde_estimate_DA_V_H}
\end{equation}

\item \label{enu:lemma:B_estimate_V_DA_H} Let $\boldsymbol{u}%
\in V,\boldsymbol{v}\in D(A),\boldsymbol{w}\in H$ then
\begin{equation}
\abs{ \left( B\left( \boldsymbol{u},\boldsymbol{v}\right) ,%
\boldsymbol{w}\right) } \leq c\vnorm{ \boldsymbol{u}
} \vnorm{ \boldsymbol{v}} ^{1/2}\lnorm{ A%
\boldsymbol{v}} ^{1/2}\lnorm{ \boldsymbol{w}} .
\label{eq:B_estimate_V_DA_H}
\end{equation}

\item \label{enu:lemma:B_estimate_DA_H_V} Let $\boldsymbol{u}\in D(A),%
\boldsymbol{v}\in H$, $\boldsymbol{w}\in V$, then
\begin{equation}
\abs{ \langle B\left( \boldsymbol{u},\boldsymbol{v}\right) ,%
\boldsymbol{w}\rangle _{V^{\prime }}} \leq c\vnorm{ \boldsymbol{u}}
^{1/2}\lnorm{ A\boldsymbol{u}}
^{1/2}\lnorm{ \boldsymbol{v}} \vnorm{ \boldsymbol{w}%
} .  \label{eq:B_estimate_DA_H_V}
\end{equation}

\item \label{enu:lemma:Btilde_estimate_V_V_V} Let $\boldsymbol{u},%
\boldsymbol{v},\boldsymbol{w}\in V$, then
\begin{equation}
\abs{ \langle \tilde{B}\left( \boldsymbol{u},\boldsymbol{v}\right) ,%
\boldsymbol{w}\rangle _{V^{\prime }}} \leq c\vnorm{ \boldsymbol{u}}
\vnorm{ \boldsymbol{v}} \lnorm{ \boldsymbol{w}} ^{1/2}\vnorm{
\boldsymbol{w}} ^{1/2}. \label{eq:Btilde_estimate_V_V_V}
\end{equation}

\item \label{enu:lemma:BandBtilde_estimate_H_V_D(A)} Let $\boldsymbol{u}\in
H$, $\boldsymbol{\boldsymbol{v}}\in V,$ $\boldsymbol{w}\in D\left(
A\right) $ and let $X$ be either $B$ or $\tilde{B}$ then
\begin{equation}
\abs{ \left\langle X\left( \boldsymbol{u},\boldsymbol{\boldsymbol{v}}%
\right) ,\boldsymbol{w}\right\rangle _{D\left( A\right) ^{\prime }}}
\leq c\lnorm{ \boldsymbol{u}} \vnorm{ \boldsymbol{\boldsymbol{v}}}
\vnorm{ \boldsymbol{w}} ^{1/2}\lnorm{ A\boldsymbol{w}} ^{1/2}.
\label{eq:BandBtilde_estimate_H_V_D(A)}
\end{equation}

\item \label{enu:lemma:Btilde_estimate_V_H_D(A)}Let $\boldsymbol{u}\in V,$ $%
\boldsymbol{\boldsymbol{v}}\in H,$ $\boldsymbol{w}\in D\left(
A\right) $ then
\begin{equation}
\abs{ \left\langle \tilde{B}\left( \boldsymbol{u},\boldsymbol{%
\boldsymbol{v}}\right) ,\boldsymbol{w}\right\rangle _{D\left(
A\right) ^{\prime }}} \leq c\left( \lnorm{ \boldsymbol{u}}
^{1/2}\vnorm{ \boldsymbol{u}} ^{1/2}\lnorm{ \boldsymbol{%
\boldsymbol{v}}} \lnorm{ A\boldsymbol{w}} +\lnorm{
\boldsymbol{\boldsymbol{v}}} \vnorm{ \boldsymbol{u}}
\vnorm{ \boldsymbol{w}} ^{1/2}\lnorm{ A\boldsymbol{w}%
} ^{1/2}\right) ,  \label{eq:Btilde_estimate_V_H_D(A)}
\end{equation}%
and hence by Poincar\'{e} inequality,
\begin{equation}
\abs{ \left\langle \tilde{B}\left( \boldsymbol{u},\boldsymbol{%
\boldsymbol{v}}\right) ,\boldsymbol{w}\right\rangle _{D\left(
A\right) ^{\prime }}} \leq c\left( \lambda _{1}\right)
^{-1/4}\vnorm{ \boldsymbol{u}} \lnorm{ \boldsymbol{\boldsymbol{v}}}
\lnorm{ A\boldsymbol{w}} . \label{eq:Btilde_estimate_V_H_D(A)_short}
\end{equation}

\item \label{enu:lemma:Btilde_estimate_D(A)_H_V}Let $\boldsymbol{u}\in
D\left( A\right) ,$ $\boldsymbol{\boldsymbol{v}}\in H,$
$\boldsymbol{w}\in V$ then
\begin{equation}
\abs{ \left\langle \tilde{B}\left( \boldsymbol{u},\boldsymbol{%
\boldsymbol{v}}\right) ,\boldsymbol{w}\right\rangle _{V^{\prime }}}
\leq c\left( \vnorm{ \boldsymbol{u}}
^{1/2}\lnorm{ A\boldsymbol{u}} ^{1/2}\lnorm{ \boldsymbol{%
\boldsymbol{v}}} \vnorm{ \boldsymbol{w}} +\lnorm{ A%
\boldsymbol{u}} \lnorm{ \boldsymbol{\boldsymbol{v}}}
\lnorm{ \boldsymbol{w}} ^{1/2}\vnorm{ \boldsymbol{w}%
} ^{1/2}\right) .  \label{eq:Btilde_estimate_D(A)_H_V}
\end{equation}
\end{enumerate}
In this lemma and throughout the paper $c$ denotes a generic scale
invariant constant.
\end{lemma}

\begin{proof}
The proof of \eqref{enu:lemma:BandBtilde_estimate_V_V_V} can be
found, for example, in \cite{b_CF88, b_T84,b_T95}  for $B$ and in
\cite[Lemma 1\textit{(iii)}]{a_FHT02} for $\tilde{B}$.

To prove \eqref{enu:lemma:BandBtilde_estimate_DA_V_H} we first
consider the case where
$\boldsymbol{u},\boldsymbol{v},\boldsymbol{w}\in \V$
\begin{align*}
&\abs{ \left( B\left( \boldsymbol{u},\boldsymbol{v}\right) ,%
\boldsymbol{w}\right) }   =\abs{
\int_{\Omega }(\boldsymbol{u}\cdot \nabla )\boldsymbol{v}\cdot \boldsymbol{w}%
dx},  \\
&| ( \tilde{B}\left( \boldsymbol{u},\boldsymbol{v}\right) ,%
\boldsymbol{w}) |  =\abs{
\int_{\Omega }\boldsymbol{u}\times (\nabla \times \boldsymbol{v})\cdot \boldsymbol{w}%
dx},
\end{align*}%
hence
\begin{align*}
\abs{ \left( X\left( \boldsymbol{u},\boldsymbol{v}\right) ,%
\boldsymbol{w}\right) } & \leq c\norm{ \boldsymbol{u}} _{L^{\infty
}}\norm{ \nabla \boldsymbol{v}} _{L^{2}}\norm{ \boldsymbol{w}}
_{L^{2}}.
\end{align*}
 By Agmon's inequality in three-dimensional space, see, e.g., \cite{b_CF88},
\begin{equation*}
\norm{ \phi } _{L^{\infty }}\leq \norm{ \phi } _{H^{1}}^{1/2}\norm{
\phi } _{H^{2}}^{1/2}
\end{equation*}%
we obtain
\begin{align*}
\abs{ \left( X\left( \boldsymbol{u},\boldsymbol{v}\right) ,%
\boldsymbol{w}\right) }
& \leq c\vnorm{ \boldsymbol{u}} ^{1/2}\lnorm{ A\boldsymbol{u}%
} ^{1/2}\vnorm{ \boldsymbol{v}} \lnorm{ \boldsymbol{w}} .
\end{align*}%
Since $\V$ is dense in $D(A)$, $V$ and $H$ we conclude the proof of %
\eqref{enu:lemma:BandBtilde_estimate_DA_V_H}.

To prove \eqref{enu:lemma:B_estimate_V_DA_H} we recall the following
Sobolev-Ladyzhenskaya inequalities (see, e.g., \cite{b_CF88,b_L85})
in 3D
\begin{align*}
\norm{ \phi } _{L^{6}}& \leq c\vnorm{ \phi
} , \\
\norm{ \phi } _{L^{3}}& \leq c\lnorm{ \phi } ^{1/2}\vnorm{ \phi }
^{1/2},
\end{align*}%
for ${\phi }\in \V$. Then we have
\begin{align*}
\abs{ \left( B\left( \boldsymbol{u},\boldsymbol{v}\right) ,%
\boldsymbol{w}\right) } & =\abs{ \int_{\Omega }(%
\boldsymbol{u}\cdot \nabla )\boldsymbol{v}\cdot \boldsymbol{w}%
dx}  \\
& \leq c\norm{ \boldsymbol{u}} _{L^{6}}\norm{
\nabla \boldsymbol{v}} _{L^{3}}\norm{ \boldsymbol{w}%
} _{L^{2}} \\
& \leq c\vnorm{ \boldsymbol{u}} \lnorm{ \nabla \boldsymbol{v}%
} ^{1/2}\vnorm{ \nabla \boldsymbol{v}}
^{1/2}\lnorm{ \boldsymbol{w}}  \\
& \leq c\vnorm{ \boldsymbol{u}} \vnorm{ \boldsymbol{v}%
} ^{1/2}\lnorm{ A\boldsymbol{v}} ^{1/2}\lnorm{ \boldsymbol{w}} .
\end{align*}
The proof of \eqref{enu:lemma:B_estimate_DA_H_V} is a direct result of the %
\eqref{enu:lemma:BandBtilde_estimate_DA_V_H} due to the symmetry
\eqref{eq:B_id1}.
The proof of \eqref{enu:lemma:Btilde_estimate_V_V_V}, %
\eqref{enu:lemma:BandBtilde_estimate_H_V_D(A)}, %
\eqref{enu:lemma:Btilde_estimate_V_H_D(A)}, %
\eqref{enu:lemma:Btilde_estimate_D(A)_H_V} can be found in
\cite[Lemma 1 \textit{(iii,iv,v,vi)}]{a_FHT02}.
\end{proof}

Using the above notations and the identity
\eqref{eq:gen_3D_vector_id} we apply $P_{\sigma }$ to
\eqref{grp:alphaMHD}  to obtain, as for the case of the NSE, the
equivalent system of equations (see, e.g., \cite{b_T84} and
\cite{a_DL72})
\begin{subequations}
\label{grp:alphaMHD:Projected}
\begin{align}
& \frac{d\boldsymbol{v}}{dt}+\tilde{B}\left( \boldsymbol{u},\boldsymbol{v}%
\right) +\nu A\boldsymbol{v}={B}\left(
\boldsymbol{B},\boldsymbol{B}\right) ,
\label{eq:alphaMHD:Projected:velocity} \\
& \frac{d\boldsymbol{B}}{dt}+B\left(
\boldsymbol{u},\boldsymbol{B}\right) -
B\left( \boldsymbol{B},\boldsymbol{u}\right) +\eta A\boldsymbol{%
B}=0,  \label{eq:alphaMHD:Projected:magField} \\
& \boldsymbol{u}(0)=\boldsymbol{u}^{in}, \\
& \boldsymbol{B}(0)=\boldsymbol{B}^{in}.
\end{align}
\end{subequations}

\begin{definition}
Let $T>0$. A weak solution of \eqref{grp:alphaMHD:Projected} in the
interval $[0,T]$, given \mbox{$\boldsymbol{u}\left( 0\right)
=\boldsymbol{u}^{in} \in V$} (or equivalently
$\boldsymbol{v}^{in}\in V^{\prime }$) and $
 \boldsymbol{B}\left( 0\right)
=\boldsymbol{B}^{in} \in H, $ is a pair of functions
$\boldsymbol{u},\ \boldsymbol{B}$, such that
\begin{equation*}
\boldsymbol{u}\in C\left( \left[ 0,T\right] ;V\right) \cap
L^{2}\left( \left[ 0,T\right] ;D\left( A\right) \right) %
\, \text{with} \,%
 \frac{d\boldsymbol{u}}{dt}\in L^{2}\left(
\left[ 0,T\right] ;H\right)
\end{equation*}
(or equivalently \mbox{$\boldsymbol{v}\in C\left( \left[ 0,T\right]
;V^{\prime
}\right) \cap L^{2}\left( \left[ 0,T\right] ;H\right) $} with \mbox{$\frac{d%
\boldsymbol{v}}{dt}\in L^{2}\left( \left[ 0,T\right] ;D\left(
A\right) ^{\prime }\right) $}) and
\begin{equation*}
 \boldsymbol{B}\in C\left( \left[ 0,T\right] ;H\right) \cap
L^{2}\left( \left[ 0,T\right] ;V\right)  %
\, \text{with} \,%
\frac{d\boldsymbol{B}}{dt}\in L^{2}\left( \left[ 0,T\right]
;V^{\prime }\right) ,
\end{equation*}
 satisfying
\begin{subequations}
\label{grp:alphaMHD:weakSol}
\begin{align}
& \left\langle \frac{d}{dt}\boldsymbol{\boldsymbol{v}},\boldsymbol{w}%
\right\rangle _{D\left( A\right)  ^{\prime} }+\left\langle
\tilde{B}\left(
\boldsymbol{u},\boldsymbol{\boldsymbol{v}}\right) ,\boldsymbol{w}%
\right\rangle _{D\left( A\right) ^{\prime} }+\nu \left( \boldsymbol{v},A%
\boldsymbol{w}\right) =\left\langle B\left( \boldsymbol{B},\boldsymbol{B}%
\right) ,\boldsymbol{w}\right\rangle _{V ^{\prime} },
\label{eq:alphaMHD:weakSol_vel} \\
& \left\langle \frac{d}{dt}\boldsymbol{B},\boldsymbol{\xi
}\right\rangle
_{V^{\prime }}+\left( B\left( \boldsymbol{u},\boldsymbol{B}\right) ,%
\boldsymbol{\xi }\right) -\left( B\left( \boldsymbol{B},\boldsymbol{u}%
\right) ,\boldsymbol{\xi }\right) +\eta \left( \left( \boldsymbol{B},%
\boldsymbol{\xi }\right) \right) =0
\label{eq:alphaMHD:weakSol_magfld}
\end{align}%
\end{subequations}
for every $\boldsymbol{w}\in D\left( A\right) ,\ \boldsymbol{\xi
}\in V$ and for almost every $t\in \left[ 0,T\right] $.

The equation \eqref{grp:alphaMHD:weakSol} is understood in the
following sense: for almost every \mbox{$t_{0},t\in \left[
0,T\right] $}
\begin{subequations}
\label{grp:alphaMHD:weakSol_integralFormulation}
\begin{align}
& \left( \boldsymbol{\boldsymbol{v}}\left( t\right)
,\boldsymbol{w}\right)
-\left( \boldsymbol{\boldsymbol{v}}\left( t_{0}\right) ,\boldsymbol{w}%
\right) +\int_{t_{0}}^{t}\left\langle \tilde{B}\left(
\boldsymbol{u}\left(
s\right) ,\boldsymbol{\boldsymbol{v}}\left( s\right) \right) ,\boldsymbol{w}%
\right\rangle _{D\left( A\right) ^{\prime }}ds+\nu
\int_{t_{0}}^{t}\left( \boldsymbol{v}\left( s\right)
,A\boldsymbol{w}\right) ds
\label{eq:alphaMHD:weakSol_integralFormulation_vel} \\
& \qquad \qquad =\int_{t_{0}}^{t}\left\langle B\left(
\boldsymbol{B}\left( s\right) ,\boldsymbol{B}\left( s\right) \right)
,\boldsymbol{w}\right\rangle
_{V^{\prime }}ds,  \notag \\
& \left( \boldsymbol{B}\left( t\right) ,\boldsymbol{\xi }\right)
-\left( \boldsymbol{B}\left( t_{0}\right) ,\boldsymbol{\xi }\right)
+\int_{t_{0}}^{t}\left( B\left( \boldsymbol{u}\left( s\right) ,\boldsymbol{B}%
\left( s\right) \right) ,\boldsymbol{\xi }\right)
ds-\int_{t_{0}}^{t}\left( B\left( \boldsymbol{B}\left( s\right)
,\boldsymbol{u}\left( s\right) \right) ,\boldsymbol{\xi }\right) ds
\label{eq:alphaMHD:xeakSol_integralFormulation_magfld} \\
& \qquad \qquad +\eta \int_{t_{0}}^{t}\left( \left(
\boldsymbol{B}\left( s\right) ,\boldsymbol{\xi }\right) \right)
ds=0.  \notag
\end{align}
\end{subequations}

When $\boldsymbol{u}^{in}\in D(A)$ (or equivalently
$\boldsymbol{v}^{in}\in H$) and $\boldsymbol{B}^{in}\in V$ we call a
strong solution of \eqref{grp:alphaMHD:Projected} in the interval
$[0,T]$ the solution that satisfies
\begin{align*}
\boldsymbol{B}\in C\left( \left[ 0,T\right] ;V\right) \cap
L^{2}\left( \left[ 0,T\right] ;D(A)\right), \,\, \boldsymbol{u}\in
C\left( \left[ 0,T\right] ;D(A)\right) \cap L^{2}( \left[ 0,T\right]
;D(A^{3/2}))
\end{align*}
(or equivalently $\boldsymbol{v}\in C\left( \left[ 0,T\right]
;H\right) \cap L^{2}\left( \left[ 0,T\right] ;D(V)\right)$).
\begin{equation*}
\end{equation*}
\end{definition}

\section{Global existence and uniqueness}

In this section we show the global well-posedness of the MHD-$\alpha
$ model \eqref{grp:alphaMHD} or equivalently
\eqref{grp:alphaMHD:Projected}.

\begin{theorem}
\label{thm:alphaMHD:weakSol} Let $\boldsymbol{u}^{in}\in V,\,\boldsymbol{B}%
^{in}\in H$. Then for any $T>0$ there exists a unique weak solution $%
\boldsymbol{u},\boldsymbol{B}$ of \eqref{grp:alphaMHD:Projected} on
$\left[ 0,T\right] $. Moreover, this solution satisfies
\begin{equation*}
\boldsymbol{u}\in L_{loc}^{\infty }\left( \left( 0,T\right]
;H^{3}\left( \Omega \right) \right),
\end{equation*}
as well as the energy equality
\begin{multline}\label{eq:alphaMHD:weakSol_energyEquality}
\lnorm{ \boldsymbol{u}\left( t\right) } ^{2}+\alpha ^{2}\vnorm{
\boldsymbol{u}\left( t\right) } ^{2}+\lnorm{ \boldsymbol{B} \left(
t\right) } ^{2}
+2 \int_{t_0}^{t}\left( \nu(\vnorm{ \boldsymbol{u}(s)} ^{2}+\alpha
^{2}\lnorm{ A\boldsymbol{u}_(s)} ^{2})+\eta\vnorm{ \boldsymbol{B}
(s)} ^{2}\right)ds
\\
= \lnorm{ \boldsymbol{u}\left( t_0\right) } ^{2}+\alpha ^{2}\vnorm{
\boldsymbol{u}\left( t_0\right) } ^{2}+\lnorm{ \boldsymbol{B}\left(
t_0\right) } ^{2}, \qquad 0\leq t_0\leq t \leq T.
\end{multline}
\end{theorem}
We use the Galerkin approximation scheme to prove the global
existence and to establish the necessary \textit{a priori}
estimates. Let $ \{w_{j}\}_{j=1}^{\infty }$ be an orthonormal basis
of $H$ consisting of eigenfunctions of the operator $A$. Denote
$H_m=\operatorname{span}\{w_1 ,\ldots, w_m\}$ and let $P_{m}$ be the
$L^{2}$-orthogonal projection from $H$ onto $H_{m}$. The Galerkin
approximation of \eqref{grp:alphaMHD:Projected} is the ordinary
differential system
\begin{subequations}
\label{grp:alphaMHD:Galerkin}
\begin{align}
& \frac{d\boldsymbol{v}_{m}}{dt}+P_{m}\tilde{B}\left( \boldsymbol{u}_{m},%
\boldsymbol{v}_{m}\right) +\nu Av_{m}=P_{m}B\left( \boldsymbol{B}_{m},%
\boldsymbol{B}_{m}\right)  \label{eq:alphaMHD:Galerkin:velocity} \\
& \frac{d\boldsymbol{B}_{m}}{dt}+P_{m}B\left( \boldsymbol{u}_{m},\boldsymbol{%
B}_{m}\right) -P_{m}B\left(
\boldsymbol{B}_{m},\boldsymbol{u}_{m}\right)
+\eta A\boldsymbol{B}_{m}=0  \label{eq:alphaMHD:Galerkin:magField} \\
& \boldsymbol{v}_{m}=\boldsymbol{u}_{m}+\alpha ^{2}A\boldsymbol{u}_{m} \\
& \boldsymbol{u}_{m}\left( 0\right) =P_{m}\boldsymbol{u}^{in} \\
& \boldsymbol{B}_{m}\left( 0\right) =P_{m}\boldsymbol{B}^{in}.
\end{align}%
\end{subequations}
Since the nonlinear terms are quadratic, hence locally Lipschitz,
then by the classical theory of ordinary differential equations,
system \eqref{grp:alphaMHD:Galerkin} has a unique solution for a
short interval of time $(-\tau _{m},T_{m})$. Our goal is to show
that the solutions of \eqref{grp:alphaMHD:Galerkin} remains finite
for all positive times, which implies that $T_{m}=\infty $.


\subsection{ $H^{1}$-Estimate of $\boldsymbol{u}_{m}$, $L^{2}$-Estimate of $
\boldsymbol{B}_{m}$}


We take the inner product of \eqref{eq:alphaMHD:Galerkin:velocity} with $%
\boldsymbol{u}_{m}$ and the inner product of %
\eqref{eq:alphaMHD:Galerkin:magField} with $\boldsymbol{B}_{m}$ and use %
\eqref{eq:B_id2},\eqref{eq:Btilda_id2},\eqref{eq:B_id1} to obtain
\begin{subequations}
\begin{align}
& \frac{1}{2}\frac{d}{dt}\left( \lnorm{ \boldsymbol{u}_{m}}
^{2}+\alpha ^{2}\vnorm{ \boldsymbol{u}_{m}} ^{2}\right) +\nu
\left( \vnorm{ \boldsymbol{u}_{m}} ^{2}+\alpha ^{2}\lnorm{ A%
\boldsymbol{u}_{m}} ^{2}\right) =\left( B\left( \boldsymbol{B}_{m},%
\boldsymbol{B}_{m}\right) ,\boldsymbol{u}_{m}\right) ,
\label{eq:alphaMHD:velocity:inner_product_u_m} \\
& \frac{1}{2}\frac{d}{dt}\lnorm{ \boldsymbol{B}_{m}} ^{2}+\eta
\vnorm{ \boldsymbol{B}_{m}} ^{2}=-\left( B\left( \boldsymbol{B}%
_{m},\boldsymbol{B}_{m}\right) ,\boldsymbol{u}_{m}\right) .
\label{eq:alphaMHD:magField:inner_product_B_m}
\end{align}%
\end{subequations}
Now, by summing up \eqref{eq:alphaMHD:velocity:inner_product_u_m} and %
\eqref{eq:alphaMHD:magField:inner_product_B_m}, we have
\begin{equation}\label{eq:alphaMHD:u+B}
\frac{1}{2}\frac{d}{dt}\left( \lnorm{ \boldsymbol{u}_{m}}
^{2}+\alpha ^{2}\vnorm{ \boldsymbol{u}_{m}} ^{2}+\lnorm{
\boldsymbol{B}_{m}} ^{2}\right) +\nu \left( \vnorm{ \boldsymbol{%
u}_{m}} ^{2}+\alpha ^{2}\lnorm{ A\boldsymbol{u}_{m}} ^{2}\right)
+\eta \vnorm{ \boldsymbol{B}_{m}} ^{2}=0.
\end{equation}%
We denote $\mu =\min \left\{ \nu ,\eta \right\} $ and obtain
\begin{equation}
\frac{1}{2}\frac{d}{dt}\left( \lnorm{ \boldsymbol{u}_{m}}
^{2}+\alpha ^{2}\vnorm{ \boldsymbol{u}_{m}} ^{2}+\lnorm{
\boldsymbol{B}_{m}} ^{2}\right) +\mu \left( \vnorm{ \boldsymbol{%
u}_{m}} ^{2}+\alpha ^{2}\lnorm{ A\boldsymbol{u}_{m}} ^{2}+\vnorm{
\boldsymbol{B}_{m}} ^{2}\right) \leq 0. \label{eq:alphaMHD:u+B
inequality}
\end{equation}%
Using Poicar\'{e}'s inequality we get
\begin{equation*}
\frac{d}{dt}\left( \lnorm{ \boldsymbol{u}_{m}} ^{2}+\alpha
^{2}\vnorm{ \boldsymbol{u}_{m}} ^{2}+\lnorm{ \boldsymbol{B}%
_{m}} ^{2}\right) +2\mu \lambda _{1}\left( \lnorm{ \boldsymbol{u%
}_{m}} ^{2}+\alpha ^{2}\vnorm{ \boldsymbol{u}_{m}} ^{2}+\lnorm{
\boldsymbol{B}_{m}} ^{2}\right) \leq 0.
\end{equation*}%
and then by Gronwall's inequality we obtain
\begin{equation*}
\lnorm{ \boldsymbol{u}_{m}\left( t\right) } ^{2}+\alpha ^{2}\vnorm{
\boldsymbol{u}_{m}\left( t\right) } ^{2}+\lnorm{
\boldsymbol{B}_{m}\left( t\right) } ^{2}\leq e^{-2\mu \lambda
_{1}t}\left( \lnorm{ \boldsymbol{u}_{m}\left( 0\right) } ^{2}+\alpha
^{2}\vnorm{ \boldsymbol{u}_{m}\left( 0\right) } ^{2}+\lnorm{
\boldsymbol{B}_{m}\left( 0\right) } ^{2}\right) .
\end{equation*}%
Hence
\begin{equation}\label{eq:alphaMHD:H1_estimate}
\lnorm{ \boldsymbol{u}_{m}\left( t\right) } ^{2}+\alpha ^{2}\vnorm{
\boldsymbol{u}_{m}\left( t\right) } ^{2}+\lnorm{
\boldsymbol{B}_{m}\left( t\right) } ^{2}\leq k_{1}:=\lnorm{
\boldsymbol{u}^{in}} ^{2}+\alpha ^{2}\vnorm{ \boldsymbol{u}%
^{in}} ^{2}+\lnorm{ \boldsymbol{B}^{in}} ^{2},
\end{equation}%
for all $t\geq 0$.

This implies that $T_m=\infty$. Indeed, consider $[0,T_m^{max})$,
the maximal interval of existence. Either $T_m^{max}=\infty$ and we
are done, or $T_m^{max}<\infty$ and we have $\lim \sup_{t\rightarrow
\left(T_m^{max}\right)^{-} } \left(\lnorm{ \boldsymbol{u}_{m}\left(
t\right) } ^{2}+\lnorm{ \boldsymbol{B}_{m}\left( t\right) }
^{2}\right)=\infty$, a contradiction to
\eqref{eq:alphaMHD:H1_estimate}. Hence we have global existence of
$\boldsymbol u_m,\,\boldsymbol B_m$, and hereafter we take an
arbitrary interval $[0,T]$.

Integrating \eqref{eq:alphaMHD:u+B} over the interval $\left(
s,t\right) $ and using the estimate \eqref{eq:alphaMHD:H1_estimate}
we obtain that, for all $0\leq s\leq t$,
\begin{equation}
2 \int_{s}^{t}\left( \nu(\vnorm{ \boldsymbol{u}_{m}(\tau)}
^{2}+\alpha
^{2}\lnorm{ A\boldsymbol{u}_{m}(\tau)} ^{2})+\eta\vnorm{ \boldsymbol{B}%
_{m}(\tau)} ^{2}\right)d\tau \leq k_{1}.
\label{eq:alphaMHD:integral_H2_estimate}
\end{equation}

\subsection{ $H^{2}$-Estimate of $\boldsymbol{u}_{m}$, $H^{1}$-Estimate of $
\boldsymbol{B}_{m}$}
By taking the inner product of \eqref{eq:alphaMHD:Galerkin:velocity}
with $A \boldsymbol{u}_{m}$ and the inner product of
\eqref{eq:alphaMHD:Galerkin:magField} with $A\boldsymbol{B}_{m}$ we
have
\begin{subequations}
\begin{align}
& \frac{1}{2}\frac{d}{dt}\left( \vnorm{ \boldsymbol{u}_{m}}
^{2}+\alpha ^{2}\lnorm{ A\boldsymbol{u}_{m}} ^{2}\right) +\nu \left(
\lnorm{ A\boldsymbol{u}_{m}} ^{2}+\alpha ^{2}\lnorm{
A^{3/2}\boldsymbol{u}_{m}} ^{2}\right) =\left( B\left(
\boldsymbol{B}_{m},\boldsymbol{B}_{m}\right)
,A\boldsymbol{u}_{m}\right) -\left( \tilde{B}\left(
\boldsymbol{u}_{m},\boldsymbol{v}_{m}\right)
,A\boldsymbol{u}_{m}\right) ,
\label{eq:alphaMHD:velocity:inner_product_Au_m} \\
& \frac{1}{2}\frac{d}{dt}\vnorm{ \boldsymbol{B}_{m}} ^{2}+\eta
\lnorm{ A\boldsymbol{B}_{m}} ^{2}=\left( B\left(
\boldsymbol{B}_{m},\boldsymbol{u}_{m}\right)
,A\boldsymbol{B}_{m}\right) -\left( B\left(
\boldsymbol{u}_{m},\boldsymbol{B}_{m}\right)
,A\boldsymbol{B}_{m}\right) .
\label{eq:alphaMHD:magField:inner_product_AB_m}
\end{align}%
\end{subequations}
First, we estimate the nonlinear terms. By
\eqref{eq:Btilde_estimate_V_V_V} we have
\begin{equation}
\abs{ \left( \tilde{B}\left( \boldsymbol{u}_{m},\boldsymbol{v}%
_{m}\right) ,A\boldsymbol{u}_{m}\right) } \leq c\left( \lambda
_{1}^{-1}+\alpha ^{2}\right) \vnorm{ \boldsymbol{u}_{m}}
\lnorm{ A\boldsymbol{u}_{m}} ^{1/2}\lnorm{ A^{3/2}%
\boldsymbol{u}_{m}} ^{3/2}.  \label{eq:Btilde_estimate_um_vm_Aum}
\end{equation}%
To bound the term $\abs{ \left( B\left( \boldsymbol{B}_{m},
\boldsymbol{B}_{m}\right) ,A\boldsymbol{u}_{m}\right) } $ we use
\eqref{eq:B_estimate_DA_H_V}
\begin{equation}
\abs{ \left( B\left( \boldsymbol{B}_{m},\boldsymbol{B}_{m}\right) ,A%
\boldsymbol{u}_{m}\right) } \leq c\vnorm{ \boldsymbol{B}%
_{m}} ^{1/2}\lnorm{ A\boldsymbol{B}_{m}}
^{1/2}\lnorm{ \boldsymbol{B}_{m}} \lnorm{ A^{3/2}\boldsymbol{%
u}_{m}} .  \label{eq:B_estimate_Bm_Bm_Aum}
\end{equation}%
By \eqref{eq:BandBtilde_estimate_DA_V_H} we have
\begin{equation}
\abs{ \left( B\left( \boldsymbol{B}_{m},\boldsymbol{u}_{m}\right) ,A
\boldsymbol{B}_{m}\right) } \leq c\vnorm{ \boldsymbol{B}%
_{m}} ^{1/2}\vnorm{ \boldsymbol{u}_{m}} \lnorm{ A
\boldsymbol{B}_{m}} ^{3/2}  \label{eq:B_estimate_Bm_um_ABm}
\end{equation}%
and by \eqref{eq:B_estimate_V_DA_H}
\begin{equation} \label{eq:B_estimate_um_Bm_ABm}
\abs{ \left( B\left( \boldsymbol{u}_{m},\boldsymbol{B}_{m}\right) ,A
\boldsymbol{B}_{m}\right) } \leq c\vnorm{ \boldsymbol{B}%
_{m}} ^{1/2}\vnorm{ \boldsymbol{u}_{m}} \lnorm{ A
\boldsymbol{B}_{m}} ^{3/2}.
\end{equation}%
Now, summing up \eqref{eq:alphaMHD:velocity:inner_product_Au_m}
 and \eqref{eq:alphaMHD:magField:inner_product_AB_m}, we obtain
\begin{multline} \label{eq:alphaMHD:sum:inner_product_A}
\frac{1}{2}\frac{d}{dt}\left( \vnorm{ \boldsymbol{u}_{m}}
^{2}+\alpha ^{2}\lnorm{ A\boldsymbol{u}_{m}} ^{2}+\vnorm{
\boldsymbol{B}_{m}} ^{2}\right) +\nu \left( |{A\boldsymbol{u}_{m}|}%
^{2}+\alpha ^{2}|{A^{3/2}\boldsymbol{u}_{m}|}^{2}\right) +\eta |{A%
\boldsymbol{B}_{m}|}^{2} \\
=\left( B\left( \boldsymbol{B}_{m},\boldsymbol{B}_{m}\right) ,A\boldsymbol{u}%
_{m}\right) -\left( \tilde{B}\left( \boldsymbol{u}_{m},\boldsymbol{v}%
_{m}\right) ,A\boldsymbol{u}_{m}\right) +\left( B\left( \boldsymbol{B}_{m},%
\boldsymbol{u}_{m}\right) ,A\boldsymbol{B}_{m}\right) -\left(
B\left( \boldsymbol{u}_{m},\boldsymbol{B}_{m}\right)
,A\boldsymbol{B}_{m}\right) .
\end{multline}%
By \eqref{eq:Btilde_estimate_um_vm_Aum}, \eqref{eq:B_estimate_Bm_Bm_Aum}, %
\eqref{eq:B_estimate_Bm_um_ABm} and \eqref{eq:B_estimate_um_Bm_ABm}
and several applications of Young's inequality we reach
\begin{multline}\label{eq:alphaMHD:u_m_inequality:H2}
\frac{d}{dt}\left( \vnorm{ \boldsymbol{u}_{m}} ^{2}+\alpha
^{2}\lnorm{ A\boldsymbol{u}_{m}} ^{2}+\vnorm{ \boldsymbol{B}%
_{m}} ^{2}\right) +\nu \left( |{A\boldsymbol{u}_{m}|}^{2}+\alpha
^{2}|{A^{3/2}\boldsymbol{u}_{m}|}^{2}\right) +\eta |{A\boldsymbol{B}_{m}|}%
^{2}  \\
\leq c(\alpha ^{2}\nu )^{-3}\left( \lambda _{1}^{-1}+\alpha
^{2}\right)
^{4}\vnorm{ \boldsymbol{u}_{m}} ^{4}|{A\boldsymbol{u}_{m}|}%
^{2}+c(\alpha ^{2}\nu )^{-2}\eta ^{-1}\vnorm{ \boldsymbol{B}%
_{m}} ^{2}\lnorm{ \boldsymbol{B}_{m}} ^{4}+c\eta
^{-3}\vnorm{ \boldsymbol{B}_{m}} ^{2}\vnorm{ \boldsymbol{u}%
_{m}} ^{4},
\end{multline}%
Integrating over $(s,t)$ and using \eqref{eq:alphaMHD:H1_estimate}, %
\eqref{eq:alphaMHD:integral_H2_estimate} we obtain
\begin{multline}\label{eq:alphaMHD:H2_estimate_of_u, with_u(s)}
\vnorm{ \boldsymbol{u}_{m}\left( t\right) } ^{2}+\alpha ^{2}\lnorm{
A\boldsymbol{u}_{m}\left( t\right) } ^{2}+\vnorm{
\boldsymbol{B}_{m}\left( t\right) } ^{2}%
+ {\int_{s}^{t}}\left( {\nu \left( |{A\boldsymbol{u}_{m}}\left( \tau \right) {|%
}^{2}+\alpha ^{2}|{A^{3/2}\boldsymbol{u}_{m}\left( \tau \right) |}%
^{2}\right) +\eta |{A\boldsymbol{B}_{m}\left( \tau \right)
|}^{2}}\right) d\tau \\
 \leq \vnorm{ \boldsymbol{u}_{m}\left( s\right)
} ^{2}+\alpha ^{2}\lnorm{ A\boldsymbol{u}_{m}\left( s\right) }
^{2}+\vnorm{ \boldsymbol{B}_{m}\left( s\right) } ^{2}+{K}_{1},
\end{multline}
where we denote
\begin{equation*}
{K}_{1} := c\left( \left( \lambda _{1}^{-1}+\alpha ^{2}\right)
^{4}\nu ^{-4}\alpha ^{-12}+\eta ^{-2}\alpha ^{-4}\left( \nu
^{-2}+\eta ^{-2}\right) \right) k_{1}^{3}.
\end{equation*}%

\begin{enumerate}
\item Now, if $\boldsymbol{u}^{in}\in D\left( A\right) $, $\boldsymbol{B}%
^{in}\in V$, we have
\begin{multline}
\vnorm{ \boldsymbol{u}_{m}\left( t\right) } ^{2}+\alpha ^{2}\lnorm{
A\boldsymbol{u}_{m}\left( t\right) }
^{2}+\vnorm{ \boldsymbol{B}_{m}\left( t\right) } ^{2}\\
+{\int_{0}^{t}}\left( {\nu \left( |{A\boldsymbol{u}_{m}}\left( \tau \right) {|%
}^{2}+\alpha ^{2}|{A^{3/2}\boldsymbol{u}_{m}\left( \tau \right) |}%
^{2}\right) +\eta |{A\boldsymbol{B}_{m}\left( \tau \right)
|}^{2}}\right) d\tau \\
\leq
\vnorm{ {\boldsymbol{u}}^{in}} ^{2}+\alpha ^{2}\lnorm{ A%
\boldsymbol{u}^{in}} ^{2}+\vnorm{ \boldsymbol{B}%
^{in}} ^{2}+{K}_{1}:={k}_{2}. \label{eq:alphaMHD:H2 and integral_H3
estimates,uin_in_D(A),Bin_in_V}
\end{multline}%

\item Otherwise, if $\boldsymbol{u}^{in}\notin D\left( A\right) $, $%
\boldsymbol{B}^{in}\notin V$, we integrate
\begin{equation*}
\vnorm{ \boldsymbol{u}_{m}\left( t\right) } ^{2}+\alpha ^{2}\lnorm{
A\boldsymbol{u}_{m}\left( t\right) } ^{2}+\vnorm{
\boldsymbol{B}_{m}\left( t\right) } ^{2}
 \leq \vnorm{ \boldsymbol{u}_{m}\left( s\right)
} ^{2}+\alpha ^{2}\lnorm{ A\boldsymbol{u}_{m}\left( s\right) }
^{2}+\vnorm{ \boldsymbol{B}_{m}\left( s\right) } ^{2}+{K}_{1}
\end{equation*}
 with respect to $s$ over $
(0,t)$ and use \eqref{eq:alphaMHD:integral_H2_estimate} to obtain
\begin{equation*}
t\left( \vnorm{ \boldsymbol{u}_{m}\left( t\right) } ^{2}+\alpha
^{2}\lnorm{ A\boldsymbol{u}_{m}\left( t\right) } ^{2}+\vnorm{
\boldsymbol{B}_{m}\left( t\right) } ^{2}\right) \leq \frac{1}{2\mu
}k_{1}+{K}_{1}t,
\end{equation*}%
hence for $t>0$

\begin{equation}
{\vnorm{ \boldsymbol{u}_{m}\left( t\right) } ^{2}+\alpha ^{2}\lnorm{
A\boldsymbol{u}_{m}\left( t\right) }
^{2}+\vnorm{ {B}_{m}\left( t\right) } ^{2}\leq {K}_{1}+%
\frac{1}{2t}\mu ^{-1}k_{1}:={k}_{2}\left( t\right) },
\label{eq:alphaMHD:H2_estimate_of_u,H1_estimate_of_B}
\end{equation}%
and thus
\begin{equation}
{\int_{s}^{t}\left( {\nu \left( |{A\boldsymbol{u}_{m}}\left( \tau \right) {|}%
^{2}+\alpha ^{2}|{A^{3/2}\boldsymbol{u}_{m}\left( \tau \right)
|}^{2}\right) +\eta |{A\boldsymbol{B}_{m}\left( \tau \right)
|}^{2}}\right) d\tau \leq 2 {K}_{1}{+\frac{1}{2t}k_{1}\mu
^{-1}}={K}_{1}+{k}_{2}\left( s\right) }.
\label{eq:alphaMHD:integral_H3_estimate}
\end{equation}
\end{enumerate}


\subsection{\protect $H^{3}$-Estimate of $\boldsymbol{u}_{m}$}
We establish a uniform upper bound for the $H^{3}$-norm of
$\boldsymbol{u}_{m}$ by providing the estimate for the vorticity
\mbox{$\boldsymbol{q}_{m}=\nabla \times \boldsymbol{v}_{m}$}. The
Galerkin approximation \eqref{eq:alphaMHD:Galerkin:velocity} is
equivalent to
\begin{equation*}
\frac{d\boldsymbol{v}_{m}}{dt}+\nu A\boldsymbol{v}_{m}-P_{m}\left(
\boldsymbol{u}_{m}\times \boldsymbol{q}_{m}\right) =P_{m}B\left( \boldsymbol{%
B}_{m},\boldsymbol{B}_{m}\right) .
\label{eq:alphaMHD:Galerkin:vorticity}
\end{equation*}%
Taking the curl of the above equation we obtain%
\begin{equation}\label{eq:alphaMHD:q:curl}
\frac{d\boldsymbol{q}_{m}}{dt}+\nu A\boldsymbol{q}_{m}-\nabla \times
P_{m}\left( \boldsymbol{u}_{m}\times \boldsymbol{q}_{m}\right)
=\nabla \times P_{m}B\left(
\boldsymbol{B}_{m},\boldsymbol{B}_{m}\right) .
\end{equation}%
We use that in periodic boundary conditions
\begin{equation}\label{eq:gen_3D_vector_id_per_bnd_cond}
\int_\Omega (\nabla \times \phi)\cdot \psi dx = \int_\Omega \phi
\cdot (\nabla \times \psi) dx
\end{equation}
and for divergence free vectors
\begin{equation}\label{eq:gen_3D_vector_div_free}
\nabla \times(\phi \times \psi)
=-(\phi\cdot\nabla)\psi+(\psi\cdot\nabla)\phi.
\end{equation}
 Taking the inner product of \eqref{eq:alphaMHD:q:curl} with
$\boldsymbol{q}_{m}$, using that $\nabla \cdot \boldsymbol{q}_{m}=0$
and the identities \eqref{eq:gen_3D_vector_id_per_bnd_cond},
\eqref{eq:gen_3D_vector_div_free} and \eqref{eq:B_id2}, we
reach%
\begin{equation*}
\frac{1}{2}\frac{d}{dt}\lnorm{ \boldsymbol{q}_{m}} ^{2}+\nu
\vnorm{ \boldsymbol{q}_{m}} ^{2}=\left( B\left( \boldsymbol{q}%
_{m},\boldsymbol{u}_{m}\right) ,\boldsymbol{q}_{m}\right) +\left(
B\left(
\boldsymbol{B}_{m},\boldsymbol{B}_{m}\right) ,\nabla \times \boldsymbol{q}%
_{m}\right) .
\end{equation*}%
We bound the right hand side using
\eqref{eq:BandBtilde_estimate_V_V_V},  Young's inequality and
\eqref{eq:alphaMHD:H1_estimate}
\begin{align*}
\abs{ \left( B\left( \boldsymbol{q}_{m},\boldsymbol{u}_{m}\right) ,%
\boldsymbol{q}_{m}\right) } & \leq c\lnorm{ \boldsymbol{q}%
_{m}} ^{1/2}\vnorm{ \boldsymbol{u}_{m}} \vnorm{
\boldsymbol{q}_{m}} ^{3/2} \\
& \leq c\nu ^{-3}\alpha ^{-4}k_{1}^{2}\lnorm{ \boldsymbol{q}%
_{m}} ^{2}+\frac{\nu }{4}\vnorm{ \boldsymbol{q}_{m}} ^{2}\text{ }
\end{align*}%
and by \eqref{eq:BandBtilde_estimate_DA_V_H}
\begin{align}\label{eq:alphaMHD:vorticity_BtildaEstimate}
\abs{ \left( B\left( \boldsymbol{B}_{m},\boldsymbol{B}_{m}\right)
,\nabla \times \boldsymbol{q}_{m}\right) } & \leq c\vnorm{
\boldsymbol{B}_{m}} ^{3/2}\lnorm{ A\boldsymbol{B}%
_{m}} ^{1/2}\vnorm{ \boldsymbol{q}_{m}} \\
& \leq c\nu ^{-1}\lnorm{ A\boldsymbol{B}_{m}} ^{2}+\vnorm{
\boldsymbol{B}_{m}} ^{6}+\frac{\nu }{4}\vnorm{ \boldsymbol{q}%
_{m}} ^{2}. \notag
\end{align}%
Note that since $\nabla \cdot \boldsymbol{v}_{m}=0$ and due to the
periodic boundary conditions we have
\begin{equation*}
\lnorm{ \boldsymbol{q}_{m}} =\lnorm{ \nabla \times
\boldsymbol{v}_{m}} =\lnorm{ \nabla \boldsymbol{v}%
_{m}} =\vnorm{ \boldsymbol{v}_{m}} ,
\end{equation*}%
hence
\begin{equation}
\lnorm{ \boldsymbol{q}_{m}} ^{2}\leq \vnorm{ \boldsymbol{%
\boldsymbol{u}}_{m}+\alpha ^{2}A\boldsymbol{\boldsymbol{u}}_{m}}
^{2}\leq \left( \lambda _{1}^{-1}+\alpha ^{2}\right) ^{2}\lnorm{ A^{3/2}%
\boldsymbol{\boldsymbol{u}}_{m}} ^{2}\text{.}
\label{eq:alphaMHD:vorticityL2norm}
\end{equation}%
Hence we obtain%
\begin{equation}
\frac{1}{2}\frac{d}{dt}\lnorm{ \boldsymbol{q}_{m}} ^{2}+\frac{%
\nu }{2}\vnorm{ \boldsymbol{q}_{m}} ^{2}\leq c\nu ^{-3}\alpha
^{-4}k_{1}^{2}\left( \lambda _{1}^{-1}+\alpha ^{2}\right)
^{2}\lnorm{
A^{3/2}\boldsymbol{\boldsymbol{u}}_{m}} ^{2}+c\nu ^{-1}\lnorm{ A%
\boldsymbol{B}_{m}} ^{2}+\vnorm{ \boldsymbol{B}_{m}} ^{6}.
\label{eq:alphaMHD:H3 u+B inequality}
\end{equation}%
In the following we denote by $c_{i}$ some constants depending on
$\nu,\eta ,\alpha ,k_{1},\lambda _{1}$. Integrating over $\left(
s,t\right) $ and
using \eqref{eq:alphaMHD:integral_H3_estimate} and %
\eqref{eq:alphaMHD:H2_estimate_of_u,H1_estimate_of_B} we have
\begin{equation}
\lnorm{ \boldsymbol{q}_{m}\left( t\right) } ^{2}\leq \lnorm{
\boldsymbol{q}_{m}\left( s\right) } ^{2}+c_{0}\left( 2{K}_{1}
+\frac{1}{2s}k_{1}\mu ^{-1}\right) +2\int_{s}^{t}\left(
{K_{1}+\frac{1 }{2\tau }k_{1}\mu ^{-1}}\right) ^{3}d\tau .
\label{eq:alphaMHD:vorticity_intermediate}
\end{equation}%
We integrate this expression with respect to $s$ over $\left( \frac{t}{2}%
,t\right) $, $t>0$ and use \eqref{eq:alphaMHD:integral_H3_estimate},
\eqref{eq:alphaMHD:vorticityL2norm} to obtain
\begin{equation}
{\lnorm{ \boldsymbol{q}_{m}\left( t\right) } ^{2}\leq \frac{1}{2%
}K_{1}^{3}t+c_{1}+\frac{c_{2}}{t}+\frac{c_{3}}{t^{2}}}
\label{eq:alphaMHD:H3_estimate_unbounded}
\end{equation}%
For $t>\frac{1}{\nu \lambda _{1}}$ we integrate %
\eqref{eq:alphaMHD:vorticity_intermediate} with respect to $s$ over
the interval $\left( t-\frac{1}{\nu \lambda _{1}},t\right) $. Note
that, by
applying also \eqref{eq:alphaMHD:integral_H3_estimate} and %
\eqref{eq:alphaMHD:vorticityL2norm}, we have
\begin{equation}
\lnorm{ \boldsymbol{q}_{m}\left( t\right) } ^{2}\leq c_{4}+c_{5}%
{\left( t-\frac{1}{\nu \lambda _{1}}\right) }^{-1}+c_{6}{\ln }\left( {1-%
\frac{1}{\nu \lambda _{1}t}}\right) ^{-1}.
\label{eq:alphaMHD:H3_estimate_bounded}
\end{equation}

From \eqref{eq:alphaMHD:H3_estimate_unbounded} and %
\eqref{eq:alphaMHD:H3_estimate_bounded} we have, for $t>0$,
\begin{equation}  \label{eq:alphaMHD:H3_estimate}
\lnorm{ \boldsymbol{q}_{m}\left( t\right) } ^{2}\leq k_{3}\left(
t\right),
\end{equation}%
where $k_{3}\left( t\right) $ has the following properties

\begin{enumerate}
\item $k_{3}\left( t\right) $ is finite for all $t>0$;

\item $k_{3}\left( t\right) $ is independent of $m$;

\item If either $\boldsymbol{\boldsymbol{u}}^{in}\notin D\left(
A^{3/2}\right) $ or $\boldsymbol{B}^{in}\notin V$, then $k_{3}\left(
t\right) $ depends on $\nu ,\eta ,\alpha ,\lnorm{ \boldsymbol{\boldsymbol{%
u}}^{in}} ,\vnorm{ \boldsymbol{\boldsymbol{u}}^{in}} ,\lnorm{
\boldsymbol{B}^{in}} $ and \mbox{$\lim_{t\rightarrow
0^{+}}k_{3}\left( t\right) =\infty $};

\item $\lim \sup_{t\rightarrow \infty }k_{3}\left( t\right) =R^{2}<\infty $,
$R^{2}$ depends on $\nu ,\eta ,\alpha $, but not on $\boldsymbol{\boldsymbol{%
u}}^{in}$ and $\boldsymbol{B}^{in}$.
\end{enumerate}

Returning to \eqref{eq:alphaMHD:H3 u+B inequality} and integrating over $%
\left(t,t+\tau\right)$, for $t>0$, $\tau\geq 0$ and using %
\eqref{eq:alphaMHD:H3_estimate} we obtain
\begin{equation}  \label{eq:alphaMHD:integral_H4_estimate}
\nu \int_t^{t+\tau} \vnorm{ \boldsymbol{q}_m} ^2 \leq k_4(t,\tau),
\end{equation}
where $k_4(t,\tau)$ as a function of $t$ satisfies properties (i)-(iii) as $%
k_3(t)$ above.

\begin{remark}
If $\boldsymbol{B}^{in}\in V$ and $\boldsymbol{u}^{in} \in D\left(
A\right) $, then by \eqref{eq:alphaMHD:H2 and integral_H3
estimates,uin_in_D(A),Bin_in_V}, Young's and Poincar\'{e}
inequalities we can bound
\eqref{eq:alphaMHD:vorticity_BtildaEstimate} by
\begin{equation*}
\abs{ \left( B\left( \boldsymbol{B}_{m},\boldsymbol{B}_{m}\right)
,\nabla \times \boldsymbol{q}_{m}\right) } \leq c\nu
^{-1}\lambda _{1}^{-1/2}k_{2}\lnorm{ A\boldsymbol{B}_{m}} ^{2}+%
\frac{\nu }{4}\vnorm{ \boldsymbol{q}_{m}} ^{2}.
\end{equation*}%
Hence we have%
\begin{equation*}
\frac{1}{2}\frac{d}{dt}\lnorm{ \boldsymbol{q}_{m}} ^{2}+\frac{%
\nu }{2}\vnorm{ \boldsymbol{q}_{m}} ^{2}\leq c\nu ^{-3}\alpha
^{-4}k_{1}^{2}\left( \lambda _{1}^{-1}+\alpha ^{2}\right)
^{2}\lnorm{ A^{3/2}\boldsymbol{\boldsymbol{u}}_{m}} ^{2}+c\nu
^{-1}\lambda _{1}^{-1/2}k_{2}\lnorm{ A\boldsymbol{B}_{m}} ^{2}
\end{equation*}%
and by integrating over $\left( 0,t\right) $ and using
\eqref{eq:alphaMHD:H2 and integral_H3
estimates,uin_in_D(A),Bin_in_V} we obtain
\begin{equation*}
\lnorm{ \boldsymbol{q}_{m}\left( t\right) } ^{2}\leq \lnorm{
\boldsymbol{q}_{m}\left( 0\right) } ^{2}+c\nu ^{-4}\alpha
^{-6}k_{1}^{2}\left( \lambda _{1}^{-1}+\alpha ^{2}\right) ^{2}\left(
{K}_{1}+{k}_{2}\right) +c\nu ^{-1}\lambda _{1}^{-1/2}k_{2}\eta
^{-1}\left( {K}_{1}+{k}_{2}\right).
\end{equation*}%
If, additionally,  $\boldsymbol{u}^{in}\in D\left( A^{3/2}\right) $,
then using \eqref{eq:alphaMHD:vorticityL2norm}, we obtain
\begin{equation}
\lnorm{ \boldsymbol{q}_{m}\left( t\right) } ^{2}\leq \left(
\lambda _{1}^{-1}+\alpha ^{2}\right) ^{2}\lnorm{ A^{3/2}\boldsymbol{u}%
^{in}} ^{2}+c\nu ^{-4}\alpha ^{-6}k_{1}^{2}\left( \lambda
_{1}^{-1}+\alpha ^{2}\right) ^{2}\left( {K}_{1}+{k}_{2}\right) +c\nu
^{-1}\lambda _{1}^{-1/2}k_{2}\eta ^{-1}\left( {K}_{1}+{k}_{2}\right)
. \label{eq:alphaMHD:H3_estimate_uin_in_D(A3/2)}
\end{equation}
\end{remark}


\subsection{Existence of weak solutions}


Let us summarize our estimates. For any $T>0$ we have

\begin{enumerate}
\item From \eqref{eq:alphaMHD:H1_estimate}%
\begin{equation}\label{eq:alphaMHD:um_Linf(0,T;V)_bound}
\norm{ \boldsymbol{u}_{m}} _{L^{\infty }\left( \left[ 0,T%
\right] ;H\right) }^{2}\leq {k_{1}}, \,%
\norm{ \boldsymbol{u}_{m}} _{L^{\infty }\left( \left[ 0,T%
\right] ;V\right) }^{2}\leq \frac{{k_{1}}}{{\alpha ^{2}}}
\,\,\text{or}\,\,%
\norm{ \boldsymbol{v}_{m}} _{L^{\infty }\left( \left[ 0,T%
\right] ;V^{\prime }\right) }^{2}\leq \frac{{k_{1}}}{{\alpha
^{2}}}\left( \lambda _{1}^{-1}+\alpha ^{2}\right) ^{2},%
\end{equation}
\begin{equation}\label{eq:alphaMHD:Bm_Linf(0,T;H)_bound}
\norm{ \boldsymbol{B}_{m}} _{L^{\infty }\left( \left[ 0,T%
\right] ;H\right) }^{2}\leq {k_{1}.}
\end{equation}

\item From \eqref{eq:alphaMHD:integral_H2_estimate} we have
\begin{equation}
\norm{ \boldsymbol{u}_{m}} _{L^{2}\left( \left[ 0,T%
\right] ;V\right) }^{2}\leq {\frac{{k_{1}}}{2{\nu }}},
\label{eq:alphaMHD:um_L2(0,T;V)_bound}
\end{equation}%
\begin{equation}
\norm{ \boldsymbol{u}_{m}} _{L^{2}\left( \left[ 0,T%
\right] ;D\left( A\right) \right) }^{2}\leq {\frac{{k_{1}}}{2{\nu
}\alpha ^{2}}}  \label{eq:alphaMHD:um_L2(0,T;D(A))_bound}
\end{equation}%
or%
\begin{equation}
\norm{ \boldsymbol{v}_{m}} _{L^{2}\left( \left[ 0,T%
\right] ;H\right) }^{2}\leq {\frac{{k_{1}}}{2{\nu }\alpha
^{2}}}\left( \lambda _{1}^{-1}+\alpha ^{2}\right) ^{2},
\label{eq:alphaMHD:vm_L2(0,T;H)_bound}
\end{equation}%
and%
\begin{equation}
\norm{ \boldsymbol{B}_{m}} _{L^{2}\left( \left[ 0,T%
\right] ;V\right) }^{2}\leq {\frac{{k_{1}}}{2{\eta }}}
\label{eq:alphaMHD:Bm_L2(0,T;V)_bound}.
\end{equation}

\item From \eqref{eq:alphaMHD:H2_estimate_of_u,H1_estimate_of_B} we have for any $\tau
\in \left( 0,T\right] $%
\begin{equation*}
\norm{ \boldsymbol{u}_{m}} _{L^{\infty }\left( \left[ \tau ,T\right]
;D\left( A\right) \right) }^{2}\leq \frac{{k}_{2}\left( \tau \right)
}{\alpha ^{2}}
\,\,\text{or}\,\,%
\norm{ \boldsymbol{v}_{m}} _{L^{\infty }\left( \left[
\tau ,T\right] ;H\right) }^{2}\leq \frac{{{k}_{2}\left( \tau \right) }%
}{{\alpha ^{2}}}\left( \lambda _{1}^{-1}+\alpha ^{2}\right) ^{2}
\end{equation*}%
and%
\begin{equation*}
\norm{ \boldsymbol{B}_{m}} _{L^{\infty }\left( \left[ \tau ,T\right]
;V\right) }^{2}\leq {{k}_{2}\left( \tau \right)} ,
\end{equation*}%
where ${k}_{2}\left( \tau\right) \rightarrow \infty $ as $\tau
\rightarrow 0^{+}$.
\end{enumerate}

Now we establish uniform estimates, in $m$, for $\frac{d\boldsymbol{u}_{m}}{%
dt}$, $\frac{d\boldsymbol{v}_{m}}{dt}$. Let us recall
\eqref{eq:alphaMHD:Galerkin:velocity}.
We have, by \eqref{eq:alphaMHD:vm_L2(0,T;H)_bound},%
\begin{equation*}
\norm{ A\boldsymbol{v}_{m}} _{L^{2}\left( \left[ 0,T%
\right] ;D\left( A\right) ^{\prime }\right) }^{2}\leq {\frac{{k_{1}}}{2{\nu }%
\alpha ^{2}}}\left( \lambda _{1}^{-1}+\alpha ^{2}\right) ^{2}.
\end{equation*}%
Also, by \eqref{eq:Btilde_estimate_V_H_D(A)_short},
\begin{equation*}
\norm{ P_{m}\tilde{B}\left( \boldsymbol{u}_{m},\boldsymbol{%
\boldsymbol{v}}_{m}\right) } _{D\left( A\right) ^{\prime }}\leq
c\left( \lambda _{1}\right) ^{-1/4}\vnorm{ \boldsymbol{u}_{m}}
\lnorm{ \boldsymbol{\boldsymbol{v}}_{m}} ,
\end{equation*}%
hence, applying \eqref{eq:alphaMHD:um_Linf(0,T;V)_bound} and %
\eqref{eq:alphaMHD:vm_L2(0,T;H)_bound},%
\begin{equation*}
\norm{ P_{m}\tilde{B}\left( \boldsymbol{u}_{m},\boldsymbol{%
\boldsymbol{v}}_{m}\right) } _{L^{2}\left( \left[ 0,T\right]
;D\left( A\right) ^{\prime }\right) }^{2}\leq
c{\frac{{k_{1}^{2}}\left( \lambda _{1}^{-1}+\alpha ^{2}\right)
^{2}}{{\alpha ^{4}}{\nu }\lambda _{1}^{1/2}}.}
\end{equation*}%
Additionally, by \eqref{eq:BandBtilde_estimate_H_V_D(A)}, we have%
\begin{equation*}
\norm{ P_{m}B\left( \boldsymbol{B}_{m},\boldsymbol{B}_{m}\right) }
_{D\left( A\right) ^{\prime }}\leq c\left( \lambda _{1}\right)
^{-1/4}\lnorm{ \boldsymbol{B}_{m}} \vnorm{ \boldsymbol{B}_{m}} ,
\end{equation*}%
therefore, using \eqref{eq:alphaMHD:Bm_Linf(0,T;H)_bound} and %
\eqref{eq:alphaMHD:Bm_L2(0,T;V)_bound}, we obtain%
\begin{equation*}
\norm{ P_{m}B\left( \boldsymbol{B}_{m},\boldsymbol{B}_{m}\right) }
_{L^{2}\left( \left[ 0,T\right] ;D\left( A\right) ^{\prime }\right)
}^{2}\leq c\frac{{k_{1}^{2}}}{{\eta }\lambda _{1}^{1/2}}.
\end{equation*}%
Consequently, by \eqref{eq:alphaMHD:Galerkin:velocity} and the above
\begin{equation}
\norm{ \frac{d\boldsymbol{\boldsymbol{v}}_{m}}{dt}} _{L^{2}\left(
\left[ 0,T\right] ;D\left( A\right) ^{\prime }\right)
}^{2}\leq c{\frac{{k_{1}^{2}}\left( \lambda _{1}^{-1}+\alpha ^{2}\right) ^{2}%
}{\alpha ^{4}{\nu }\lambda _{1}^{1/2}}}+ {\frac{{k_{1}}\left(
\lambda _{1}^{-1}+\alpha ^{2}\right) ^{2}}{2\alpha ^{2}}}+c\frac{{%
k_{1}^{2}}}{{\eta }\lambda _{1}^{1/2}}:=K
\label{eq:alphaMHD:dvm_dt_L2(0,T;D(A)')_bound}
\end{equation}%
and, in particular,

\begin{equation}
\norm{ \frac{d\boldsymbol{\boldsymbol{u}}_{m}}{dt}} _{L^{2}\left(
\left[ 0,T\right] ;H\right) }^{2}\leq \frac{K}{\alpha ^{4}}.
\label{eq:alphaMHD:dum_dt_L2(0,T;H)_bound}
\end{equation}

Now we establish uniform estimates, in $m$, for $\frac{d\boldsymbol{B}_{m}}{%
dt}$. Let us recall \eqref{eq:alphaMHD:Galerkin:magField}.
We have, by \eqref{eq:alphaMHD:Bm_L2(0,T;V)_bound},%
\begin{equation*}
\norm{ A\boldsymbol{B}_{m}} _{L^{2}\left( \left[ 0,T%
\right] ;V^{\prime }\right) }^{2}\leq \frac{{k_{1}}}{2{\eta}}.
\end{equation*}%
Also, by \eqref{eq:BandBtilde_estimate_V_V_V},
\begin{equation*}
\norm{ P_{m}B\left( \boldsymbol{u}_{m},\boldsymbol{B}_{m}\right) }
_{V^{\prime }}\leq c\left( \lambda _{1}\right)
^{-1/4}\vnorm{ \boldsymbol{u}_{m}} \vnorm{ \boldsymbol{B}%
_{m}} ,
\end{equation*}%
Hence, by \eqref{eq:alphaMHD:um_Linf(0,T;V)_bound} and
\eqref{eq:alphaMHD:Bm_L2(0,T;V)_bound},
\begin{align*}
\norm{ P_{m}B\left( \boldsymbol{u}_{m},\boldsymbol{B}_{m}\right) }
_{L^{2}\left( \left[ 0,T\right] ;V^{\prime }\right) }^{2}&\leq
c\frac{{k_{1}^{2}}}{2{\alpha ^{2}\eta }\lambda _{1}^{1/2}}.
\end{align*}%
Similarly%
\begin{equation*}
\norm{ P_{m}B\left( \boldsymbol{B}_{m},\boldsymbol{u}_{m}\right) }
_{V^{\prime }}\leq c\left( \lambda _{1}\right)
^{-1/4}\vnorm{ \boldsymbol{B}_{m}} \vnorm{ \boldsymbol{u}%
_{m}}
\end{equation*}%
and
\begin{equation*}
\norm{ P_{m}B\left( \boldsymbol{B}_{m},\boldsymbol{u}_{m}\right) }
_{L^{2}\left( \left[ 0,T\right] ;V^{\prime }\right) }^{2}\leq
c\frac{{k_{1}^{2}}}{2{\alpha ^{2}\eta}\lambda _{1}^{1/2}}.
\end{equation*}%
Hence, from the above and \eqref{eq:alphaMHD:Galerkin:magField}, we
have
\begin{equation}
\norm{ \frac{d\boldsymbol{B}_{m}}{dt}} _{L^{2}\left( %
\left[ 0,T\right] ;V^{\prime }\right) }^{2}\leq
c\frac{{k_{1}^{2}}}{\alpha ^{2}{\eta }\lambda _{1}^{1/2}}+
\frac{{k_{1}}}{2}:=\tilde{K}.
\label{eq:alphaMHD:dBm_dt_L2(0,T;V')_bound}
\end{equation}

From \eqref{eq:alphaMHD:um_L2(0,T;D(A))_bound} and %
\eqref{eq:alphaMHD:dum_dt_L2(0,T;H)_bound}, using Aubin's
Compactness Lemma (see, for example, \cite[Lemma
8.4]{b_CF88},\cite{b_L69} or \cite{b_T84}), we may assume that there
exists a subsequence $\boldsymbol{u}_{m^{\prime }}$ of
$\boldsymbol{u}_{m}$ and $\boldsymbol{u}\in L^{2}\left( \left[
0,T\right] ;D\left( A\right) \right) \cap C\left( \left[ 0,T\right]
;H\right) $ such that
\begin{subequations}
\begin{align}
\boldsymbol{u}_{m^{\prime }}& \rightarrow \boldsymbol{u\qquad
}\text{weakly in }L^{2}\left( \left[ 0,T\right] ;D\left( A\right)
\right),
\label{eq:alphaMHD:um_weakConvL2_D(A)} \\
\boldsymbol{u}_{m^{\prime }}& \rightarrow \boldsymbol{u}\qquad \text{%
strongly in }L^{2}\left( \left[ 0,T\right] ;V\right) \text{ and }
\label{eq:alphaMHD:um_strongConvL2_V} \\
\boldsymbol{u}_{m^{\prime }}& \rightarrow \boldsymbol{u}\qquad \text{%
strongly in }C\left( \left[ 0,T\right] ;H\right) ,
\label{eq:alphaMHD:um_strongConvC_H}
\end{align}%
\end{subequations}
as $m^{\prime }\rightarrow \infty $.
Moreover,\mbox{$\left({d}/{dt}\right) \boldsymbol{u}_{m^{\prime
}}\rightarrow \left({d}/{dt}\right)\boldsymbol{u}$} weakly
in $L^{2}\left( \left[ 0,T\right] ;H\right) $.  Or equivalently, by %
\eqref{eq:alphaMHD:vm_L2(0,T;H)_bound} and %
\eqref{eq:alphaMHD:dvm_dt_L2(0,T;D(A)')_bound}, there exists a subsequence $%
\boldsymbol{v}_{m^{\prime }}$ of $\boldsymbol{v}_{m}$ such that
\begin{subequations}
\begin{align}
\boldsymbol{v}_{m^{\prime }}& \rightarrow \boldsymbol{v\qquad
}\text{weakly in }L^{2}\left( \left[ 0,T\right] ;H\right) ,
\label{eq:alphaMHD:vm_weakConvL2_H} \\
\boldsymbol{v}_{m^{\prime }}& \rightarrow \boldsymbol{v\qquad }\text{%
strongly in }L^{2}\left( \left[ 0,T\right] ;V^{\prime }\right) ,
\label{eq:alphaMHD:vm_strongConvL2_V'} \\
\boldsymbol{v}_{m^{\prime }}& \rightarrow \boldsymbol{v\qquad }\text{%
strongly in }C\left( \left[ 0,T\right] ;D\left( A\right) ^{\prime
}\right) , \label{eq:alphaMHD:vm_strongConvC_D(A)'}
\end{align}%
\end{subequations}
\mbox{$\left({d}/{dt}\right)\boldsymbol{v}_{m^{\prime }}\rightarrow \left({d}/{dt}\right)%
\boldsymbol{v}$} weakly in $L^{2}\left( \left[ 0,T\right]
;D(A)^{\prime }\right) $, as $m^{\prime }\rightarrow \infty $, where
$\boldsymbol{v=u}+\alpha ^{2}A\boldsymbol{u}$ is in $L^{2}\left(
\left[ 0,T\right] ;H\right) \cap C\left( \left[ 0,T\right] ;D\left(
A\right)
^{\prime }\right) $. Also, by \eqref{eq:alphaMHD:Bm_L2(0,T;V)_bound} and %
\eqref{eq:alphaMHD:dBm_dt_L2(0,T;V')_bound}, there exists a subsequence $%
\boldsymbol{B}_{m^{\prime }}$ of $\boldsymbol{B}_{m}$ and
$\boldsymbol{B}\in L^{2}\left( \left[ 0,T\right] ;V\right) \cap
C\left( \left[ 0,T\right] ;V^{\prime }\right) $ such that
\begin{subequations}
\begin{align}
\boldsymbol{B}_{m^{\prime }}& \rightarrow \boldsymbol{B\qquad
}\text{weakly in }L^{2}\left( \left[ 0,T\right] ;V\right) ,
\label{eq:alphaMHD:Bm_weakConvL2_V} \\
\boldsymbol{B}_{m^{\prime }}& \rightarrow \boldsymbol{B\qquad }\text{%
strongly in }L^{2}\left( \left[ 0,T\right] ;H\right) ,
\label{eq:alphaMHD:Bm_strongConvL2_H} \\
\boldsymbol{B}_{m^{\prime }}& \rightarrow \boldsymbol{B\qquad }\text{%
strongly in }C\left( \left[ 0,T\right] ;V^{\prime }\right)
\label{eq:alphaMHD:Bm_strongConvC_V'}
\end{align}%
\end{subequations}
and \mbox{$\left({d}/{dt}\right)\boldsymbol{B}_{m^{\prime
}}\rightarrow \left({d}/{dt}\right)\boldsymbol{B}$} weakly in
$L^{2}\left( \left[ 0,T\right] ;V^{\prime }\right) $, as $m^{\prime
}\rightarrow \infty $.

Since $\boldsymbol{v}_{m^{\prime }}\rightarrow \boldsymbol{v}$ weakly in $%
L^{2}\left( \left[ 0,T\right] ;H\right) $ and strongly in
$L^{2}\left( \left[ 0,T\right] ;V^{\prime }\right) $ and
$\boldsymbol{B}_{m^{\prime }}\rightarrow \boldsymbol{B}$ weakly in
$L^{2}\left( \left[ 0,T\right] ;V\right) $ and strongly in
$L^{2}\left( \left[ 0,T\right] ;H\right) $, then there exists a set
$E\subset \left[ 0,T\right] $ of Lebesgue measure zero
and a subsequence of $\boldsymbol{v}_{m^{\prime }}$, $\boldsymbol{B}%
_{m^{\prime }}$, which we relabel $\boldsymbol{v}_{m}$,
$\boldsymbol{B}_{m}$ respectively, such that
$\boldsymbol{v}_{m}\left( s\right) \rightarrow \boldsymbol{v}\left(
s\right) $ weakly in $H$ and strongly in $V^{\prime }$
for every $s\in \left[ 0,T\right] \backslash E$, and $\boldsymbol{B}%
_{m}\left( s\right) \rightarrow \boldsymbol{B}\left( s\right) $
weakly in $V$ and strongly in $H$ for every $s\in \left[ 0,T\right]
\backslash E$.

Let $\boldsymbol{w}\in D\left( A\right) $, $\boldsymbol{\xi }\in V$,
then by taking the inner product of
\eqref{eq:alphaMHD:Galerkin:velocity} with $ \boldsymbol{w}$, and of
\eqref{eq:alphaMHD:Galerkin:magField} with $ \boldsymbol{\xi }$ and
integrating over the interval $\left[ t_{0},t\right] ,\,t,t_{0}\in
\left[ 0,T\right] $, we have
\begin{subequations}
\label{grp:alphaMHD:weakFormulation_m}
\begin{align}
\left( \boldsymbol{\boldsymbol{v}}_{_{m}}\left( t\right) ,\boldsymbol{w}%
\right) -\left( \boldsymbol{\boldsymbol{v}}_{_{m}}\left( t_{0}\right) ,%
\boldsymbol{w}\right) & +\int_{t_{0}}^{t}\left( \tilde{B}\left( \boldsymbol{u%
}_{_{m}}\left( s\right) ,\boldsymbol{\boldsymbol{v}}_{_{m}}\left(
s\right) \right) ,P_{m}\boldsymbol{w}\right) ds
\label{eq:alphaMHD:weakFormulation_m_vel} \\
& +\nu \int_{t_{0}}^{t}\left( \boldsymbol{v}_{_{m}}\left( s\right) ,A%
\boldsymbol{w}\right) ds=\int_{t_{0}}^{t}\left({B}\left(
\boldsymbol{B
}_{_{m}}\left( s\right) ,\boldsymbol{B}_{_{m}}\left( s\right) \right) ,P_{m}%
\boldsymbol{w}\right) ds,  \notag \\
\left( \boldsymbol{B}_{_{m}}\left( t\right) ,\boldsymbol{\xi
}\right) -\left( \boldsymbol{B}_{_{m}}\left( t_{0}\right)
,\boldsymbol{\xi }\right) &
+\int_{t_{0}}^{t}\left( B\left( \boldsymbol{u}_{_{m}}\left( s\right) ,%
\boldsymbol{B}_{_{m}}\left( s\right) \right) ,P_{m}\boldsymbol{\xi
}\right)
ds  \label{eq:alphaMHD:weakFormulation_m_magfld} \\
& -\int_{t_{0}}^{t}\left( B\left( \boldsymbol{B}_{_{m}}\left( s\right) ,%
\boldsymbol{u}_{_{m}}\left( s\right) \right) ,P_{m}\boldsymbol{\xi
}\right)
ds+\eta \int_{t_{0}}^{t}\left( \left( \boldsymbol{B}_{_{m}}\left( s\right) ,%
\boldsymbol{\xi }\right) \right) ds=0.  \notag
\end{align}%
\end{subequations}
First we consider \eqref{eq:alphaMHD:weakFormulation_m_vel}. Since $%
\boldsymbol{v}_{m}\left( s\right) \rightarrow \boldsymbol{v}\left(
s\right) $ weakly in $H,$ then for $t,t_{0}\in \left[ 0,T\right]
\backslash E$
\begin{equation*}
\left( \boldsymbol{\boldsymbol{v}}_{_{m}}\left( t\right) ,\boldsymbol{w}%
\right) -\left( \boldsymbol{\boldsymbol{v}}_{_{m}}\left( t_{0}\right) ,%
\boldsymbol{w}\right) \rightarrow \left(
\boldsymbol{\boldsymbol{v}}\left( t\right) ,\boldsymbol{w}\right)
-\left( \boldsymbol{\boldsymbol{v}}\left( t_{0}\right)
,\boldsymbol{w}\right) ,\text{ as }m\rightarrow \infty
\end{equation*}%
and since $\boldsymbol{w}\in D\left( A\right) $  we also have
\begin{equation*}
\lim_{m\rightarrow \infty }\int_{t_{0}}^{t}\left(
\boldsymbol{v}_{m}\left(
s\right) ,A\boldsymbol{w}\right) ds=\int_{t_{0}}^{t}\left( \boldsymbol{v}%
\left( s\right) ,A\boldsymbol{w}\right) ds.
\end{equation*}%
Now
\begin{equation}
\lim_{m\rightarrow \infty }\lnorm{ P_{m}A\boldsymbol{w}-A\boldsymbol{w}%
} =\lim_{m\rightarrow \infty }\vnorm{ P_{m}\boldsymbol{w}-%
\boldsymbol{w}} =\lim_{m\rightarrow \infty }\lnorm{ P_{m}%
\boldsymbol{w}-\boldsymbol{w}} =0.  \label{eq:alphaMHD:limPmw-w}
\end{equation}%
For the nonlinear terms we have%
\begin{align*}
& \abs{ \int_{t_{0}}^{t}\left( \tilde{B}\left( \boldsymbol{u}%
_{_{m}}\left( s\right) ,\boldsymbol{\boldsymbol{v}}_{_{m}}\left(
s\right) \right) ,P_{m}\boldsymbol{w}\right) -\left\langle
\tilde{B}\left( \boldsymbol{u}\left( s\right)
,\boldsymbol{\boldsymbol{v}}\left( s\right) \right)
,\boldsymbol{w}\right\rangle _{D\left( A\right) ^{\prime
}}ds}  \\
& \quad \leq \abs{ \int_{t_{0}}^{t}\left\langle \tilde{B}\left(
\boldsymbol{u}_{_{m}}\left( s\right) ,\boldsymbol{\boldsymbol{v}}%
_{_{m}}\left( s\right) \right) ,P_{m}\boldsymbol{w}-\boldsymbol{w}%
\right\rangle _{D\left( A\right) ^{\prime }}ds}  \\
& \quad +\abs{ \int_{t_{0}}^{t}\left\langle \tilde{B}\left(
\boldsymbol{u}_{_{m}}\left( s\right) -\boldsymbol{u}\left( s\right) ,%
\boldsymbol{\boldsymbol{v}}_{_{m}}\left( s\right) \right) ,\boldsymbol{w}%
\right\rangle _{D\left( A\right) ^{\prime }}ds}  \\
& \quad +\abs{ \int_{t_{0}}^{t}\left\langle \tilde{B}\left(
\boldsymbol{u}\left( s\right)
,\boldsymbol{\boldsymbol{v}}_{_{m}}\left(
s\right) -\boldsymbol{\boldsymbol{v}}\left( s\right) \right) ,\boldsymbol{w}%
\right\rangle _{D\left( A\right) ^{\prime }}ds}  \\
& \quad =:I_{m}^{(1)}+I_{m}^{(2)}+I_{m}^{(3)}
\end{align*}%
By \eqref{eq:Btilde_estimate_V_H_D(A)_short}%
\begin{equation*}
I_{m}^{(1)}\leq c\left( \lambda _{1}\right)
^{-1/4}\int_{t_{0}}^{t}\vnorm{ \boldsymbol{u}_{_{m}}\left( s\right)
} \lnorm{ \boldsymbol{\boldsymbol{v}}_{_{m}}\left( s\right) }
\lnorm{ P_{m}A\boldsymbol{w}-A\boldsymbol{w}} ds,
\end{equation*}%
using Cauchy-Schwarz inequality we obtain%
\begin{equation*}
I_{m}^{(1)}  \leq c\left( \lambda _{1}\right) ^{-1/4}\lnorm{ P_{m}A%
\boldsymbol{w}-A\boldsymbol{w}} \norm{ \boldsymbol{u}%
_{m}} _{L^{2}\left( \left[ 0,T\right] ;V\right) }\norm{
\boldsymbol{v}_{m}} _{L^{2}\left( \left[ 0,T\right] ;H\right) },
\end{equation*}%
hence by \eqref{eq:alphaMHD:um_L2(0,T;V)_bound},
\eqref{eq:alphaMHD:vm_L2(0,T;H)_bound} and
\eqref{eq:alphaMHD:limPmw-w} $ \lim_{m\rightarrow \infty
}I_{m}^{(1)}=0$.

Again, by \eqref{eq:Btilde_estimate_V_H_D(A)_short},
\begin{equation*}
I_{m}^{(2)}\leq c\left( \lambda _{1}\right)
^{-1/4}\int_{t_{0}}^{t}\vnorm{ \boldsymbol{u}_{m}\left( s\right) -
\boldsymbol{u}\left( s\right) } \lnorm{ \boldsymbol{v}_{m}\left(
s\right) } \lnorm{ A\boldsymbol{w}} ds,
\end{equation*}%
and by Cauchy-Schwarz and \eqref{eq:alphaMHD:vm_L2(0,T;H)_bound},
\begin{equation*}
I_{m}^{(2)}\leq c\left( \lambda _{1}\right) ^{-1/4}\lnorm{
A\boldsymbol{w} } {\frac{{k_{1}}}{2{\nu }\alpha ^{2}}}\left( \lambda
_{1}^{-1}+\alpha ^{2}\right) ^{2}\left( \int_{0}^{T}\vnorm{
\boldsymbol{u}_{m}\left( s\right) -\boldsymbol{u}\left( s\right) }
^{2}ds\right) ^{1/2},
\end{equation*}%
hence $\lim_{m\rightarrow \infty }I_{m}^{(2)}=0$, since $
\boldsymbol{u}_{m}\rightarrow \boldsymbol{u}$ in $L^{2}\left( \left[
0,T \right] ;V\right) $.

Finally, we show that $\lim_{m\rightarrow \infty }I_{m}^{(3)}=0$. We
define a linear functional for $\boldsymbol{h }\in L^{2}\left(
\left[ 0,T\right]
;H\right) $ by%
\begin{equation*}
\phi\left( \boldsymbol{\boldsymbol{h}}\right) ={\
\int_{t_{0}}^{t}\left\langle \tilde{B}\left( \boldsymbol{u}\left( s\right) ,%
\boldsymbol{\boldsymbol{h}}\right) ,\boldsymbol{w}\right\rangle
_{D\left( A\right) ^{\prime }}ds} ,
\end{equation*}%
by \eqref{eq:Btilde_estimate_V_H_D(A)_short} and Cauchy-Schwarz
\begin{equation*}
\abs{ \phi \left( \boldsymbol{\boldsymbol{h}}\right) }
\leq c\left( \lambda _{1}\right) ^{-1/4}\lnorm{ A\boldsymbol{w}%
} \vnorm{ \boldsymbol{u}\left( s\right) }
_{L^{2}\left( \left[ 0,T\right] ;V\right) }\norm{ \boldsymbol{%
\boldsymbol{h}}\left( s\right) } _{L^{2}\left( \left[ 0,T\right]
;H\right) }
\end{equation*}%
hence, due to \eqref{eq:alphaMHD:um_L2(0,T;V)_bound}, $\phi $ is a
bounded linear functional, and thus, since
$\boldsymbol{v}_{m}\rightarrow \boldsymbol{v}$ weakly in
$L^{2}\left( \left[ 0,T\right] ;H\right) $,
\begin{equation*}
\lim_{m\rightarrow \infty }{\phi \left( \boldsymbol{v}_{m}\left( s\right) -%
\boldsymbol{v}\left( s\right) \right) }=0.
\end{equation*}%
and hence $\lim_{m\rightarrow \infty }I_{m}^{(3)}=0$. It remains to
pass to the limit in the right hand side element of
\eqref{eq:alphaMHD:weakFormulation_m_vel}.
\begin{align*}
& \abs{ \int_{t_{0}}^{t}\left( B\left( \boldsymbol{B}_{_{m}}\left(
s\right) ,\boldsymbol{B}_{_{m}}\left( s\right) \right) ,P_{m}\boldsymbol{w}%
\right) -\left\langle B\left( \boldsymbol{B}\left( s\right) ,\boldsymbol{B}%
\left( s\right) \right) ,\boldsymbol{w}\right\rangle _{V^{\prime
}}ds}  \\
& \quad \leq \abs{ \int_{t_{0}}^{t}\left\langle {B}\left(
\boldsymbol{B}_{_{m}}\left( s\right) ,\boldsymbol{B}_{_{m}}\left(
s\right) \right) ,P_{m}\boldsymbol{w}-\boldsymbol{w}\right\rangle
_{V^{\prime
}}ds}  \\
& \quad +\abs{ \int_{t_{0}}^{t}\left\langle {B}\left(
\boldsymbol{B}_{_{m}}\left( s\right) -\boldsymbol{B}\left( s\right) ,%
\boldsymbol{B}_{_{m}}\left( s\right) \right)
,\boldsymbol{w}\right\rangle
_{V^{\prime }}ds}  \\
& \quad +\abs{ \int_{t_{0}}^{t}\left\langle {B}\left(
\boldsymbol{B}\left( s\right) ,\boldsymbol{B}_{_{m}}\left( s\right) -%
\boldsymbol{B}\left( s\right) \right) ,\boldsymbol{w}\right\rangle
_{V^{\prime }}ds}  \\
& \quad =:J_{m}^{(1)}+J_{m}^{(2)}+J_{m}^{(3)}.
\end{align*}%
Now, by \eqref{eq:BandBtilde_estimate_V_V_V} and Poincar\'{e}
inequality
\begin{equation*}
J_{m}^{(1)}\leq c\left( \lambda _{1}\right) ^{-1/4}\vnorm{ P_{m}%
\boldsymbol{w}-\boldsymbol{w}} \norm{ \boldsymbol{B}%
_{_{m}}} _{L^{2}\left( \left[ 0,T\right] ;V\right) }^{2},
\end{equation*}%
hence, by \eqref{eq:alphaMHD:limPmw-w}, $\lim_{m\rightarrow \infty
}J_{m}^{(1)}=0$.

By \eqref{eq:BandBtilde_estimate_V_V_V} and Poincar\'{e} inequality
\begin{equation*}
J_{m}^{(2)}\leq c\left( \lambda _{1}\right)
^{-1/4}\int_{t_{0}}^{t}\vnorm{ \boldsymbol{B}_{_{m}}\left( s\right) -%
\boldsymbol{B}\left( s\right) } \vnorm{ \boldsymbol{B}%
_{_{m}}\left( s\right) } \vnorm{ \boldsymbol{w}} ds
\end{equation*}%
and, applying Cauchy-Schwarz and
\eqref{eq:alphaMHD:Bm_L2(0,T;V)_bound},
\begin{equation*}
J_{m}^{(2)}\leq c\left( \lambda _{1}\right) ^{-1/4}\vnorm{ \boldsymbol{w}%
} \frac{k_{1}}{2\eta }\left( \int_{0}^{T}\vnorm{ \boldsymbol{B}%
_{_{m}}\left( s\right) -\boldsymbol{B}\left( s\right) }
^{2}ds\right) ^{1/2},\quad
\end{equation*}%
hence, since $\boldsymbol{B}_{m}\rightarrow \boldsymbol{B}$ weakly
in $ L^{2}\left( \left[ 0,T\right] ;V\right) $, we have
$\lim_{m\rightarrow \infty }J_{m}^{(2)}=0$ (similarly to the
argument given for $I_{m}^{(3)}$).

Similarly we can show that $\lim_{m\rightarrow \infty
}J_{m}^{(3)}=0$.

It remains to pass to the limit in
\eqref{eq:alphaMHD:weakFormulation_m_magfld}. Note
that%
\begin{equation}
\lim_{m\rightarrow \infty }\lnorm{ P_{m}\boldsymbol{\xi
}-\boldsymbol{\xi }} =0.  \label{eq:alphaMHD:limPmx-x}
\end{equation}%
We recall that $\boldsymbol{B}_{m}\left( s\right) \rightarrow \boldsymbol{B}%
\left( s\right) $ weakly in $V$ and strongly in $H$ for every $s\in
\left[ 0,T\right] \backslash E$, hence the convergence for the
linear terms is easy. For the nonlinear terms we have
\begin{align*}
& \abs{ \int_{t_{0}}^{t}\left( B\left( \boldsymbol{u}_{_{m}}\left(
s\right) ,\boldsymbol{B}_{_{m}}\left( s\right) \right)
,P_{m}\boldsymbol{\xi }\right) -\left( B\left( \boldsymbol{u}\left(
s\right) ,\boldsymbol{B}\left(
s\right) \right) ,\boldsymbol{\xi }\right) ds}  \\
& \quad \leq \abs{ \int_{t_{0}}^{t}\left( B\left( \boldsymbol{u}%
_{_{m}}\left( s\right) ,\boldsymbol{B}_{_{m}}\left( s\right) \right) ,P_{m}%
\boldsymbol{\xi }-\boldsymbol{\xi }\right) ds}  \\
& \quad +\abs{ \int_{t_{0}}^{t}\left( B\left( \boldsymbol{u}%
_{_{m}}\left( s\right) -\boldsymbol{u}\left( s\right) ,\boldsymbol{B}%
_{_{m}}\left( s\right) \right) ,\boldsymbol{\xi }\right) ds}  \\
& \quad +\abs{ \int_{t_{0}}^{t}\left( B\left( \boldsymbol{u}\left(
s\right) ,\boldsymbol{B}_{_{m}}\left( s\right) -\boldsymbol{B}\left(
s\right) \right) ,\boldsymbol{\xi }\right) ds}  \\
& \quad =:S_{m}^{(1)}+S_{m}^{(2)}+S_{m}^{(3)}.
\end{align*}%
Now, by \eqref{eq:BandBtilde_estimate_DA_V_H} and Cauchy-Schwarz
inequality
\begin{equation*}
S_{m}^{(1)}\leq c\left( \lambda _{1}\right) ^{-1/4}\lnorm{ P_{m}%
\boldsymbol{\xi }-\boldsymbol{\xi }} \norm{ \boldsymbol{u}%
_{_{m}}} _{L^{2}\left( \left[ 0,T\right] ;D\left( A\right)
\right) }\norm{ \boldsymbol{B}_{_{m}}} _{L^{2}\left( %
\left[ 0,T\right] ;V\right) },
\end{equation*}%
hence, by \eqref{eq:alphaMHD:um_L2(0,T;D(A))_bound},%
\eqref{eq:alphaMHD:Bm_L2(0,T;V)_bound} and \eqref{eq:alphaMHD:limPmx-x}, $%
\lim_{m\rightarrow \infty }S_{m}^{(1)}=0$.

Again, by \eqref{eq:BandBtilde_estimate_DA_V_H} and Cauchy-Schwarz
inequality

\begin{equation*}
S_{m}^{(2)}\leq c\left( \lambda _{1}\right) ^{-1/4}\lnorm{ \boldsymbol{%
\xi }} \norm{ \boldsymbol{B}_{_{m}}}
_{L^{2}\left( \left[ 0,T\right] ;V\right) }\norm{ \boldsymbol{u}%
_{_{m}}\left( s\right) -\boldsymbol{u}\left( s\right) }
_{L^{2}\left( \left[ 0,T\right] ;D\left( A\right) \right) }^{2}
\end{equation*}%
hence, since $\boldsymbol{u}_{m}\rightarrow \boldsymbol{u}$ weakly in $%
L^{2}\left( \left[ 0,T\right] ;D\left( A\right) \right) $, we have  $%
\lim_{m\rightarrow \infty }S_{m}^{(2)}=0$ (similarly to the case for $%
I_{m}^{(3)}$). By similar arguments, using that $\boldsymbol{B}%
_{m}\rightarrow \boldsymbol{B}$ weakly in $L^{2}\left( \left[
0,T\right] ;V\right) $, we obtain also that $\lim_{m\rightarrow
\infty }S_{m}^{(3)}=0.$

For the term $\int_{t_{0}}^{t}\left( B\left(
\boldsymbol{B}_{_{m}}\left(
s\right) ,\boldsymbol{u}_{_{m}}\left( s\right) \right) ,P_{m}\boldsymbol{\xi}%
\right) $ we can perform the same estimates using the %
\eqref{eq:B_estimate_V_DA_H} to bound operator $B$.

Hence, we can pass to the limit in
\eqref{grp:alphaMHD:weakFormulation_m} and we obtain that for every
$t,t_{0}\in \left[ 0,T\right] \backslash E$

\begin{subequations}
\label{grp:alphaMHD:weakFormulation}
\begin{align}
\left( \boldsymbol{\boldsymbol{v}}\left( t\right)
,\boldsymbol{w}\right)
-\left( \boldsymbol{\boldsymbol{v}}\left( t_{0}\right) ,\boldsymbol{w}%
\right) +\int_{t_{0}}^{t}\left\langle \tilde{B}\left(
\boldsymbol{u}\left(
s\right) ,\boldsymbol{\boldsymbol{v}}\left( s\right) \right) ,\boldsymbol{w}%
\right\rangle _{D\left( A\right) ^{\prime }}ds& +\nu
\int_{t_{0}}^{t}\left( \boldsymbol{v}\left( s\right)
,A\boldsymbol{w}\right) ds
\label{eq:alphaMHD:weakFormulation_vel} \\
& =\int_{t_{0}}^{t}\left\langle B\left( \boldsymbol{B}\left( s\right) ,%
\boldsymbol{B}\left( s\right) \right) ,\boldsymbol{w}\right\rangle
_{V^{\prime }}ds,  \notag \\
\left( \boldsymbol{B}\left( t\right) ,\boldsymbol{\xi }\right)
-\left( \boldsymbol{B}\left( t_{0}\right) ,\boldsymbol{\xi }\right)
+\int_{t_{0}}^{t}\left( B\left( \boldsymbol{u}\left( s\right) ,\boldsymbol{B}%
\left( s\right) \right) ,\boldsymbol{\xi }\right) ds&
-\int_{t_{0}}^{t}\left( B\left( \boldsymbol{B}\left( s\right) ,\boldsymbol{u}%
\left( s\right) \right) ,\boldsymbol{\xi }\right) ds
\label{eq:alphaMHD:weakFormulation_magfld} \\
& +\eta \int_{t_{0}}^{t}\left( \left( \boldsymbol{B}\left( s\right) ,%
\boldsymbol{\xi }\right) \right) ds=0.  \notag
\end{align}%
\end{subequations}
for every $\boldsymbol{w}\in D\left( A\right) ,\boldsymbol{\xi }\in
V$.

Now we show that $\boldsymbol{\boldsymbol{v}}\in C\left( \left[
0,T\right]
;V^{\prime }\right) $ (or equivalently $\boldsymbol{u}\in C\left( \left[ 0,T%
\right] ;V\right) $) and $\boldsymbol{B}\in C\left( \left[
0,T\right] ;H\right) $.

Notice that since $\norm{ \boldsymbol{v}_{m}}
_{L^{\infty }\left( \left[ 0,T\right] ;V^{\prime }\right) }^{2}\leq \frac{{%
k_{1}}}{{\alpha ^{2}}}\left( \lambda _{1}^{-1}+\alpha ^{2}\right) ^{2}$ and $%
\boldsymbol{v}_{m}\rightarrow \boldsymbol{v}$ strongly in
$L^{2}\left( \left[
0,T\right] ;V^{\prime }\right) $ then $\norm{ \boldsymbol{v}%
} _{L^{\infty }\left( \left[ 0,T\right] ;V^{\prime }\right)
}^{2}\leq \frac{{k_{1}}}{{\alpha ^{2}}}\left( \lambda
_{1}^{-1}+\alpha ^{2}\right) ^{2}$. Hence
\eqref{eq:alphaMHD:weakFormulation_vel} implies that
$\boldsymbol{\boldsymbol{v}}\left( t\right) \in $ $C_{w}\left(
\left[
0,T\right] ;V^{\prime }\right) $ because $D\left( A\right) $ is dense in $V$%
. Since, also, for a fixed $t_{0}$, $\norm{ \boldsymbol{\boldsymbol{v}%
}\left( t\right) } _{V^{\prime }}\rightarrow \norm{
\boldsymbol{\boldsymbol{v}}\left( t_{0}\right) } _{V^{\prime }}$%
, as $t\rightarrow t_{0}$, then we have
$\boldsymbol{\boldsymbol{v}}\in
$ $C\left( \left[ 0,T\right] ;V^{\prime }\right) $, or equivalently $%
\boldsymbol{u}\in C\left( \left[ 0,T\right] ;V\right) $.

Similarly, since $\norm{ \boldsymbol{B}_{m}} _{L^{\infty
}\left( \left[ 0,T\right] ;H\right) }^{2}\leq {k_{1}}$ and $\boldsymbol{B}%
_{m}\rightarrow \boldsymbol{B}$ strongly in $L^{2}\left( \left[
0,T\right]
;H\right) $ and $D\left( A\right) $ is dense in $H$ and because of %
\eqref{eq:alphaMHD:weakFormulation_magfld} we have
$\boldsymbol{B}\in C\left( \left[ 0,T\right] ;H\right) $.


\subsection{Uniqueness and continuous dependence of weak solutions on the initial data}


Next, we show the continuous dependence of weak solutions on the
initial data and, in particular, the uniqueness of weak solutions.

Let $\boldsymbol{u},\,\boldsymbol{B}$ and $\boldsymbol{\bar{u}},\,%
\boldsymbol{\bar{B}}$ be any two weak solutions of %
\eqref{grp:alphaMHD:Projected} on the interval $\left[ 0,T\right] $
with initial values \mbox{$\boldsymbol{u}\left( 0\right)
=\boldsymbol{u}^{in}$}, \mbox{$\boldsymbol{B}\left( 0\right)
=\boldsymbol{B}^{in}$},
\mbox{$\boldsymbol{\bar{u}}%
\left( 0\right) =\boldsymbol{\bar{u}}^{in}$},
\mbox{$\boldsymbol{\bar{B}}\left( 0\right)
=\boldsymbol{\bar{B}}^{in}$}. We denote
\mbox{$\boldsymbol{v=u}+\alpha
^{2}A\boldsymbol{u}$}, \mbox{$\boldsymbol{\bar{v}=\bar{u}}+\alpha ^{2}A\boldsymbol{%
\bar{u}}$}, \mbox{$\delta \boldsymbol{u=u}-\boldsymbol{\bar{u}}$},
\mbox{$\delta
\boldsymbol{v=v}-\boldsymbol{\bar{v}}$} and \mbox{$\delta \boldsymbol{B=B}-%
\boldsymbol{\bar{B}}$}. Then \eqref{grp:alphaMHD:Projected} implies
\begin{align}
& \frac{d}{dt}\delta \boldsymbol{v}+\nu A\delta \boldsymbol{v}+\tilde{B}%
\left( \delta \boldsymbol{u},\boldsymbol{\boldsymbol{v}}\right) +\tilde{B}%
\left( \boldsymbol{\bar{u}},\delta \boldsymbol{v}\right) =B\left(
\delta \boldsymbol{B},\boldsymbol{B}\right) +B\left(
\boldsymbol{\bar{B}},\delta
\boldsymbol{B}\right) ,  \label{eq:alphaMHD:delta_v} \\
& \frac{d}{dt}\delta \boldsymbol{B}+\eta A\delta
\boldsymbol{B}=-B\left(
\delta \boldsymbol{u},\boldsymbol{B}\right) -B\left( \boldsymbol{\bar{u}}%
,\delta \boldsymbol{B}\right) +B\left( \delta \boldsymbol{B},\boldsymbol{u}%
\right) +B\left( \boldsymbol{\bar{B}},\delta \boldsymbol{u}\right) ,
\label{eq:alphaMHD:delta_B} \\
& \delta \boldsymbol{u}\left( 0\right) =\delta \boldsymbol{u}^{in}=%
\boldsymbol{u}^{in}-\boldsymbol{\bar{u}}^{in}, \\
& \delta \boldsymbol{B}\left( 0\right) =\delta \boldsymbol{B}^{in}=%
\boldsymbol{B}^{in}-\boldsymbol{\bar{B}}^{in}.
\end{align}

Since ${d\boldsymbol{v}}/{dt}\in L^{2}\left( \left[ 0,T%
\right] ;D\left( A\right) ^{\prime }\right) ,\ \delta \boldsymbol{u }\in%
L^{2}\left( \left[ 0,T\right] ;D\left( A\right) \right) $ and $d
\boldsymbol{B}/{dt}\in L^{2}\left( \left[ 0,T\right] ;V^{\prime
}\right) $, $\boldsymbol{B},\ \boldsymbol{\bar{B}},\ \delta \boldsymbol{B}%
\in L^{2}\left( \left[ 0,T\right] ;V\right) $ and due to the
identities \eqref{eq:B_id2} and \eqref{eq:Btilda_id2}, we have for
almost every $t\in \left[ 0,T\right] $
\begin{align*}
& \left\langle \frac{d}{dt}\delta \boldsymbol{v},\delta \boldsymbol{u}%
\right\rangle _{D\left( A\right) ^{\prime }}+\nu \left( \vnorm{
\delta
\boldsymbol{u}} ^{2}+\alpha ^{2}\lnorm{ A\delta \boldsymbol{u}%
} ^{2}\right) +\left\langle \tilde{B}\left( \boldsymbol{\bar{u}}%
,\delta \boldsymbol{v}\right) ,\delta \boldsymbol{u}\right\rangle
_{D\left(
A\right) ^{\prime }} \\
& \qquad \qquad \qquad \qquad =\left\langle B\left( \delta \boldsymbol{B},%
\boldsymbol{B}\right) ,\delta \boldsymbol{u}\right\rangle _{D\left(
A\right)
^{\prime }}+\left\langle B\left( \boldsymbol{\bar{B}},\delta \boldsymbol{B}%
\right) ,\delta \boldsymbol{u}\right\rangle _{V^{\prime }}, \\
& \left\langle \frac{d}{dt}\delta \boldsymbol{B},\delta \boldsymbol{B}%
\right\rangle _{V^{\prime }}+\eta \vnorm{ \delta \boldsymbol{B}%
} ^{2}=-\left\langle B\left( \delta \boldsymbol{u},\boldsymbol{B}%
\right),\delta \boldsymbol{B}\right\rangle _{V^{\prime
}}+\left\langle B\left( \delta \boldsymbol{B},\boldsymbol{u}\right)
,\delta \boldsymbol{B}\right\rangle _{V^{\prime }}+\left\langle
B\left( \boldsymbol{\bar{B}},\delta \boldsymbol{u}\right),%
\delta \boldsymbol{B}\right\rangle _{V^{\prime }}.
\end{align*}

Notice that by theorem of interpolation by Lions and Magenes, see, e.g., \cite[%
Chap. III, Lemma 1.2]{b_T84},
\begin{equation*}
\left\langle \frac{d}{dt}\delta \boldsymbol{v},\delta \boldsymbol{u}%
\right\rangle _{D\left( A\right) ^{\prime }}=\frac{d}{dt}\left(
\lnorm{ \delta \boldsymbol{u}} ^{2}+\alpha ^{2}\vnorm{ \delta
\boldsymbol{u}} ^{2}\right)
\end{equation*}%
and%
\begin{equation*}
\left\langle \frac{d}{dt}\delta \boldsymbol{B},\delta \boldsymbol{B}%
\right\rangle _{V^{\prime }}=\frac{d}{dt}\lnorm{ \delta \boldsymbol{B}%
} ^{2},
\end{equation*}%
thus we have
\begin{subequations}
\begin{align}
& \frac{d}{dt}\left( \lnorm{ \delta \boldsymbol{u}} ^{2}+\alpha
^{2}\vnorm{ \delta \boldsymbol{u}} ^{2}\right) +\nu \left( \vnorm{
\delta \boldsymbol{u}} ^{2}+\alpha ^{2}\lnorm{ A\delta
\boldsymbol{u}} ^{2}\right) +\left\langle \tilde{B}\left(
\boldsymbol{\bar{u}},\delta \boldsymbol{v}\right) ,\delta \boldsymbol{u}%
\right\rangle _{D\left( A\right) ^{\prime }}
\label{eq:alphaMHD:delta_inner_product_delta_u} \\
& \qquad \qquad \qquad \qquad =\left\langle B\left( \delta \boldsymbol{B},%
\boldsymbol{B}\right) ,\delta \boldsymbol{u}\right\rangle _{D\left(
A\right)
^{\prime }}+\left\langle B\left( \boldsymbol{\bar{B}},\delta \boldsymbol{B}%
\right) ,\delta \boldsymbol{u}\right\rangle _{V^{\prime }},  \notag \\
& \frac{d}{dt}\lnorm{ \delta \boldsymbol{B}} ^{2}+\eta \vnorm{
\delta \boldsymbol{B}} ^{2}=-\left\langle B\left(
\delta \boldsymbol{u},\boldsymbol{B}\right) ,\delta \boldsymbol{%
B}\right\rangle _{V^{\prime }}+\left\langle B\left( \delta \boldsymbol{B},%
\boldsymbol{u}\right) ,\delta \boldsymbol{B}\right\rangle
_{V^{\prime }}+\left\langle B\left( \boldsymbol{\bar{B}},\delta \boldsymbol{u%
}\right),\delta \boldsymbol{B}\right\rangle _{V^{\prime }}.
\label{eq:alphaMHD:delta_inner_product_delta_B}
\end{align}%
\end{subequations}
By summation of \eqref{eq:alphaMHD:delta_inner_product_delta_u} and %
\eqref{eq:alphaMHD:delta_inner_product_delta_B} we obtain
\begin{multline*}
\frac{d}{dt}\left( \lnorm{ \delta \boldsymbol{u}} ^{2}+\alpha
^{2}\vnorm{ \delta \boldsymbol{u}} ^{2}+\lnorm{ \delta
\boldsymbol{B}} ^{2}\right) +\nu \left( \vnorm{ \delta
\boldsymbol{u}} ^{2}+\alpha ^{2}\lnorm{ A\delta \boldsymbol{u}%
} ^{2}\right) +\eta \vnorm{ \delta \boldsymbol{B}}
^{2} \\
=-\left\langle \tilde{B}\left( \boldsymbol{\bar{u}},\delta \boldsymbol{v}%
\right) ,\delta \boldsymbol{u}\right\rangle _{D\left( A\right)
^{\prime }}+\left\langle \tilde{B}\left( \delta
\boldsymbol{B},\boldsymbol{B}\right) ,\delta
\boldsymbol{u}\right\rangle _{D\left( A\right) ^{\prime
}}+\left\langle B\left( \delta \boldsymbol{B},\boldsymbol{u}\right)
,\delta \boldsymbol{B}\right\rangle _{V^{\prime }}.
\end{multline*}%
By \eqref{eq:Btilde_estimate_D(A)_H_V} we get
\begin{equation*}
\abs{ \left\langle \tilde{B}\left( \boldsymbol{\bar{u}},\delta
\boldsymbol{v}\right) ,\delta \boldsymbol{u}\right\rangle
_{V^{\prime }}} \leq c\left( \vnorm{ \boldsymbol{\bar{u}}}
^{1/2}\lnorm{ A\boldsymbol{\bar{u}}} ^{1/2}\lnorm{ \delta
\boldsymbol{v}} \vnorm{ \delta \boldsymbol{u}}
+\lnorm{ A\boldsymbol{\bar{u}}} \lnorm{ \delta \boldsymbol{v}%
} \lnorm{ \delta \boldsymbol{u}} ^{1/2}\vnorm{ \delta
\boldsymbol{u}} ^{1/2}\right) ,
\end{equation*}%
and by applying Young's inequality
\begin{equation*}
\abs{ \left\langle \tilde{B}\left( \boldsymbol{\bar{u}},\delta
\boldsymbol{v}\right) ,\delta \boldsymbol{u}\right\rangle
_{V^{\prime
}}} \leq \frac{c}{\nu \lambda _{1}^{1/2}}\lnorm{ A%
\boldsymbol{\bar{u}}} ^{2}\left( \lnorm{ \delta \boldsymbol{u}%
} ^{2}+\alpha ^{2}\vnorm{ \delta \boldsymbol{u}}
^{2}\right) +\frac{\nu }{2}\vnorm{ \delta \boldsymbol{u}} ^{2}+%
\frac{\nu }{4}\alpha ^{2}\lnorm{ A\delta \boldsymbol{u}} ^{2}.
\end{equation*}%
By \eqref{eq:BandBtilde_estimate_H_V_D(A)} and Young's inequality we
have
\begin{equation*}
\abs{ \left\langle \tilde{B}\left( \delta \boldsymbol{B},\boldsymbol{%
B}\right) ,\delta \boldsymbol{u}\right\rangle _{D\left( A\right)
^{\prime }}} \leq c\lnorm{ \delta \boldsymbol{B}}
^{2}\vnorm{ \boldsymbol{B}} ^{2}+\frac{1}{\nu \alpha ^{2}}%
\vnorm{ \delta \boldsymbol{u}} ^{2}+\frac{\nu }{4}\alpha ^{2}\lnorm{
A\delta \boldsymbol{u}} .
\end{equation*}%
Also, by \eqref{eq:B_estimate_V_DA_H} and Young's inequality we obtain%
\begin{equation*}
\abs{ \left( B\left( \delta \boldsymbol{B},\boldsymbol{u}\right)
,\delta \boldsymbol{B}\right) } \leq \frac{c}{\eta \lambda
_{1}^{1/2}}\lnorm{ A\boldsymbol{u}} ^{2}\lnorm{ \delta
\boldsymbol{B}} ^{2}+\frac{\eta }{2}\vnorm{ \delta \boldsymbol{B}}
^{2}
\end{equation*}%
Summing up we have%
\begin{align*}
& \frac{d}{dt}\left( \lnorm{ \delta \boldsymbol{u}} ^{2}+\alpha
^{2}\vnorm{ \delta \boldsymbol{u}} ^{2}+\lnorm{ \delta
\boldsymbol{B}} ^{2}\right) +\frac{\nu }{2}\left( \vnorm{ \delta
\boldsymbol{u}} ^{2}+\alpha ^{2}\lnorm{ A\delta \boldsymbol{u}}
^{2}\right)+\frac{\eta }{2}\vnorm{ \delta \boldsymbol{B}}
^{2}  \\
& \qquad \qquad \qquad \leq \left( \frac{c}{\nu \lambda _{1}^{1/2}}
\lnorm{ A\boldsymbol{\bar{u}}} ^{2} +\frac{1}{\nu \alpha ^{4}}+c\vnorm{ \boldsymbol{B}%
} ^{2}+\frac{c}{\eta \lambda _{1}^{1/2}}\lnorm{ A\boldsymbol{u}%
} ^{2}\right) \left( \lnorm{ \delta \boldsymbol{u}} ^{2}+\alpha
^{2}\vnorm{ \delta \boldsymbol{u}} ^{2}+\lnorm{ \delta
\boldsymbol{B}} ^{2}\right) .
\end{align*}
We denote
\begin{equation*}
z\left( s\right) = \frac{c}{\nu \lambda _{1}^{1/2}}
\lnorm{ A\boldsymbol{\bar{u}}} ^{2} +\frac{1}{\nu \alpha ^{4}}+c\vnorm{ \boldsymbol{B}%
} ^{2}+\frac{c}{\eta \lambda _{1}^{1/2}}\lnorm{ A\boldsymbol{u}%
} ^{2}
\end{equation*}
and use Gronwall's
inequality to obtain%
\begin{equation}
\lnorm{ \delta \boldsymbol{u}\left( t\right) } ^{2}+\alpha
^{2}\vnorm{ \delta \boldsymbol{u}\left( t\right) } ^{2}+\lnorm{
\delta \boldsymbol{B}\left( t\right) } ^{2}\leq \left( \lnorm{
\delta \boldsymbol{u}\left( 0\right) } ^{2}+\alpha ^{2}\vnorm{
\delta \boldsymbol{u}\left( 0\right) } ^{2}+\lnorm{ \delta
\boldsymbol{B}\left( 0\right) } ^{2}\right) \exp \left(
\int_{0}^{t}z\left( s\right) ds\right) ,
\label{eq:alphaMHD:cont_depend}
\end{equation}%
since $\boldsymbol{u},\boldsymbol{\bar{u}}\in L^{2}\left( \left[
0,T\right]
;D\left( A\right) \right) $ and $\boldsymbol{B}\in L^{2}\left( \left[ 0,T%
\right] ;V\right) $ the integral $\left( \int_{0}^{t}z\left(
s\right) ds\right) $ is finite. Hence
\eqref{eq:alphaMHD:cont_depend} implies the
continuous dependence of the weak solutions of %
\eqref{grp:alphaMHD:Projected} on the initial data in any bounded
interval of time $\left[ 0,T\right] $. In particular, the solutions
are unique.


\subsection{Strong solutions}


\begin{theorem}
Let $T>0$, $\boldsymbol{u}^{in}\in V,\,\boldsymbol{B}^{in}\in H$.
Then there
exists a unique solution $\boldsymbol{u},\boldsymbol{B}$ of %
\eqref{grp:alphaMHD:Projected} on $\left[ 0,T\right] $ satisfying%
\begin{equation}
\boldsymbol{u}\in L_{loc}^{\infty }\left( \left( 0,T\right] ;D(A^{{3/2}%
})\right) \cap L_{loc}^{2}\left( \left( 0,T\right] ;D(A^{2})\right)
\cap
C\left( \left[ 0,T\right] ;V\right) \cap L^{2}\left( \left[ 0,T%
\right] ;D(A)\right)   \label{eq:alphaMHD:strong_sol_u}
\end{equation}%
and%
\begin{equation}
\boldsymbol{B}\in L_{loc}^{\infty }\left( \left( 0,T\right]
;V\right) \cap
L_{loc}^{2}\left( \left( 0,T\right] ;D(A)\right) \cap C\left( %
\left[ 0,T\right] ;H\right) \cap L^{2}\left( \left[ 0,T\right]
;V\right) . \label{eq:alphaMHD:strong_sol_B}
\end{equation}
If $\boldsymbol{B}^{in}\in V$  and $\boldsymbol{u}^{in}\in D(A)$
then the solution is the strong solution
\begin{align*}
\boldsymbol{u}&\in C\left( \left[ 0,T\right] ;D(A)\right) \cap
L^{2}( \left[ 0,T\right] ;D(A^{3/2})) ,
\\
\boldsymbol{B}&\in C\left( \left[ 0,T\right] ;V\right) \cap
L^{2}\left( \left[ 0,T\right] ;D(A)\right).
\end{align*}
If, additionally, $\boldsymbol{u}^{in}\in D(A^{3/2})$ then
\begin{equation*}
\boldsymbol{u}\in C( \left[ 0,T\right] ;D(A^{3/2})) \cap L^{2}\left(
\left[ 0,T\right] ;D(A^{2})\right) .
\end{equation*}
\end{theorem}
\begin{remark}
Following the techniques presented in \cite{a_FT89} (see also
\cite{a_FT98}) we can show that for any $t>0$ the solution is
analytic in time with values in a Gevrey class of regularity of
spatial analytic functions. As a result, we have an exponentially
fast convergence in the wave number $m$, as $m\to\infty$, in a
certain sense, of the Galerkin approximation to the unique strong
solution of \eqref{grp:alphaMHD:Projected}, see, for instance,
\cite{a_DT93,a_JMT95}. This Gevrey regularity result also implies
the exponential decay of large wavenumber modes in the dissipation
range of turbulent flows \cite{a_DT95}.

\end{remark}
\begin{proof}
We use the Galerkin estimates derived in the previous subsections
and similar ideas and compactness theorems in the corresponding
spaces to converge to the strong solution. For
\eqref{eq:alphaMHD:strong_sol_u} and
\eqref{eq:alphaMHD:strong_sol_B} we need the estimates
\eqref{eq:alphaMHD:H3_estimate}, \eqref{eq:alphaMHD:integral_H4_estimate}, %
\eqref{eq:alphaMHD:H1_estimate}, \eqref{eq:alphaMHD:integral_H2_estimate} %
and %
\eqref{eq:alphaMHD:H2_estimate_of_u,H1_estimate_of_B},
\eqref{eq:alphaMHD:integral_H3_estimate},
\eqref{eq:alphaMHD:H1_estimate},
\eqref{eq:alphaMHD:integral_H2_estimate}. For
$\boldsymbol{B}^{in}\in V, \ \boldsymbol{u}^{in}\in D(A)$ we use the
estimate \eqref{eq:alphaMHD:H2 and integral_H3
estimates,uin_in_D(A),Bin_in_V} and if   $ \boldsymbol{u}^{in}\in
D(A^{3/2})$ we use \eqref{eq:alphaMHD:H3_estimate_uin_in_D(A3/2)}.
Also, since the strong solutions are weak, by uniqueness of weak
solutions the strong solutions are unique.
\end{proof}

\section{Convergence to the solutions of MHD equations as  $\alpha
\rightarrow 0^{+}$}

We emphasize again that our point of view is that the alpha model is
to be considered as a regularizing numerical scheme. The next
theorem shows that using the \textit{a priori} estimates established
previously, one can extract subsequences of the weak solutions of
system \eqref{grp:alphaMHD:Projected}, which converge, as $\alpha
\rightarrow 0^{+}$, (in the appropriate sense defined in the
theorem) to a Leray-Hopf weak solution of the three-dimensional MHD
equations on any time interval $\left[ 0,T \right] $. For the
definition and existence of weak solutions of the
 3D MHD equations, see, for instance, \cite{a_DL72}
and \cite{a_ST83}. The notion of a Leray-Hopf weak solution of MHD
that satisfies the energy inequality \eqref{eq:MHD:energyIneq} is
inspired from  a Leray-Hopf solution of NSE and formulated in the
theorem. Also, if the initial data is smooth we prove that a
subsequence of the strong solutions of the MHD-$\alpha$ equations
converges to the unique strong solution of the 3D MHD on an interval
$\left[0,T_*(u^{in}, B^{in})\right]$, which is the interval of
existence of the strong solution.

\begin{theorem}
Let $T>0$, $\boldsymbol{u}^{in}\in V,\,\boldsymbol{B}^{in}\in H$ and
denote
by $\boldsymbol{u}_{\alpha },\,\boldsymbol{B}_{\alpha }$ and \mbox{$%
\boldsymbol{v}_{\alpha }=\boldsymbol{u}_{\alpha }+\alpha ^{2}A\boldsymbol{u}%
_{\alpha }$} the weak solution of \eqref{grp:alphaMHD:Projected} on $%
\left[ 0,T\right] $. Then there are subsequences
$\boldsymbol{u}_{\alpha _{j}},\,\boldsymbol{v}_{\alpha
_{j}},\,\boldsymbol{B}_{\alpha _{j}}$ and a pair of functions
$\boldsymbol{v},\boldsymbol{B}\in L^{\infty }\left( \left[
0,T\right] ;H\right) \cap L^{2}\left( \left[ 0,T\right] ;V\right) $
such that, as \mbox{$\alpha _{j}\rightarrow 0^{+}$},

\begin{enumerate}
\item $\boldsymbol{u}_{\alpha _{j}}\rightarrow \boldsymbol{v}$ and $%
\boldsymbol{B}_{\alpha _{j}}\rightarrow \boldsymbol{B}$ weakly in $%
L^{2}\left( \left[ 0,T\right] ;V\right) $ and strongly in
$L^{2}\left( \left[ 0,T\right] ;H\right) $,

\item $\boldsymbol{v}_{\alpha _{j}}\rightarrow \boldsymbol{v}$ weakly in $%
L^{2}\left( \left[ 0,T\right] ;H\right) $ and strongly in
$L^{2}\left( \left[ 0,T\right] ;V^{\prime }\right) $ and

\item $\boldsymbol{u}_{\alpha _{j}}\left( t\right) \rightarrow \boldsymbol{v}\left( t\right)
$ and $\boldsymbol{B}_{\alpha _{j}}\left( t\right) \rightarrow
\boldsymbol{B}\left( t\right) $ weakly in $H$ and uniformly on
$\left[ 0,T\right] $.
\end{enumerate}

Furthermore, the pair $\boldsymbol{v},\boldsymbol{B}$ is a
Leray-Hopf weak solution of the MHD equations
\begin{align*}
& \frac{d\boldsymbol{v}}{dt}+\tilde{B}\left( \boldsymbol{v},\boldsymbol{v}%
\right) +\nu A\boldsymbol{v}=B\left(
\boldsymbol{B},\boldsymbol{B}\right) ,
\\
& \frac{d\boldsymbol{B}}{dt}+B\left(
\boldsymbol{v},\boldsymbol{B}\right) -B\left(
\boldsymbol{B},\boldsymbol{v}\right) +\eta A\boldsymbol{B}=0
\end{align*}%
 with initial data $\boldsymbol{v}%
\left( 0\right) =\boldsymbol{u}^{in},\,\boldsymbol{B}\left( 0\right) =%
\boldsymbol{B}^{in}$, which satisfies the energy inequality
\begin{equation}\label{eq:MHD:energyIneq}
\lnorm{ \boldsymbol{v}\left( t\right) } ^{2}+\lnorm{ \boldsymbol{B}
\left( t\right) } ^{2} +2 \int_{t_0}^{t}\left( \nu\vnorm{
\boldsymbol{v}(s)} ^{2}+\eta\vnorm{ \boldsymbol{B} (s)}
^{2}\right)ds \leq \lnorm{ \boldsymbol{v}\left( t_0\right) }
^{2}+\lnorm{ \boldsymbol{B}\left( t_0\right) } ^{2}
\end{equation}
for almost every $t_0$, $0\leq t_0 \leq T$  and all
$t\in\left[t_0,T\right]$.
\end{theorem}

\begin{proof}
From estimates \eqref{eq:alphaMHD:H1_estimate} and %
\eqref{eq:alphaMHD:integral_H2_estimate}, by passing to the limit
(using the proof of Theorem \ref{thm:alphaMHD:weakSol}), we have
that the solution of \eqref{grp:alphaMHD:Projected} satisfies
\begin{equation*}
{\lnorm{ \boldsymbol{u}_{\alpha }\left( t\right) } ^{2}+\alpha
^{2}\vnorm{ \boldsymbol{u}_{\alpha }\left( t\right) } ^{2}+\lnorm{
\boldsymbol{B}_{\alpha }\left( t\right) } ^{2}\leq k_{1}}
\end{equation*}%
and
\begin{equation*}
2 \int_{0}^{T}\left( \nu(\vnorm{ \boldsymbol{u}_{\alpha }(t)}
^{2}+\alpha ^{2}\lnorm{ A\boldsymbol{u}_{\alpha }(t)}
^{2})+\eta\vnorm{ \boldsymbol{B}_{\alpha } (t)} ^{2}\right)dt \leq
k_{1},
\end{equation*}%
notice that since $\alpha \rightarrow 0^{+}$ we can assume that
$0<\alpha \leq L$; consequently, we can bound the right hand side by
$\tilde{k}_{1}:=\lnorm{ \boldsymbol{u}^{in}} ^{2}+L^{2}\vnorm{
\boldsymbol{u}^{in}} ^{2}+\lnorm{ \boldsymbol{B}^{in}} ^{2}$, which
is independent of $\alpha $, therefore we can extract subsequences
$\boldsymbol{u}_{\alpha _{j}},\,\boldsymbol{v}_{\alpha
_{j}},\,\boldsymbol{B}_{\alpha _{j}}$, such that
\begin{align*}
{\boldsymbol{u}_{\alpha _{j}}}\rightarrow \boldsymbol{u\quad }&
\text{weakly
in }L^{2}\left( \left[ 0,T\right] ;V\right) , \\
{\boldsymbol{v}_{\alpha _{j}}}\rightarrow \boldsymbol{v\quad }&
\text{weakly
in }L^{2}\left( \left[ 0,T\right] ;H\right) \text{ and} \\
{\boldsymbol{B}_{\alpha _{j}}}\rightarrow \boldsymbol{B\quad }&
\text{weakly in }L^{2}\left( \left[ 0,T\right] ;V\right) ,
\end{align*}%
as $\alpha _{j}\rightarrow 0^{+}$.

Now we establish uniform estimates, independent of $\alpha $, for
${d\boldsymbol{B}
_{\alpha }}/{dt}$ and ${d\boldsymbol{u}_{\alpha }}/{dt}$. From %
\eqref{eq:alphaMHD:Projected:magField} we have
\begin{equation*}
\norm{ A^{-1}\frac{d\boldsymbol{B}_{\alpha }}{dt}} \leq \norm{
A^{-1}B\left( \boldsymbol{u}_{\alpha },\boldsymbol{B}_{\alpha
}\right) } +\norm{ A^{-1}\boldsymbol{B}\left( \boldsymbol{B}_{\alpha
},\boldsymbol{u}_{\alpha }\right) } +\eta \norm{
\boldsymbol{B}_{\alpha }} ,
\end{equation*}%
notice that by \eqref{eq:BandBtilde_estimate_H_V_D(A)}%
\begin{equation*}
\lnorm{ A^{-1}B\left( \boldsymbol{u}_{\alpha
},\boldsymbol{B}_{\alpha
}\right) } \leq c\lambda _{1}^{-1/4}\lnorm{ \boldsymbol{u}%
_{\alpha }} \vnorm{ \boldsymbol{\boldsymbol{B}}_{\alpha }} ,
\end{equation*}%
hence%
\begin{align*}
\norm{ B\left( \boldsymbol{u}_{\alpha },\boldsymbol{B}_{\alpha
}\right) } _{L^{2}\left( \left[ 0,T\right] ;D\left( A\right)
^{\prime }\right) }& \leq c\lambda _{1}^{-1/2}{\int_{0}^{T}}\lnorm{
\boldsymbol{u}_{\alpha }\left( t\right) } ^{2}\vnorm{
\boldsymbol{\boldsymbol{B}}_{\alpha }\left( t\right) } ^{2}dt \\
& \leq c\lambda _{1}^{-1/2}\tilde{k}_{1}^{2} \eta^{-1}
\end{align*}%
and similarly
\begin{equation*}
\norm{ B\left( \boldsymbol{B}_{\alpha },\boldsymbol{u}_{\alpha
}\right) } _{L^{2}\left( \left[ 0,T\right] ;D\left( A\right)
^{\prime }\right) }\leq c\lambda _{1}^{-1/2}\tilde{k}_{1}^{2}\eta
^{-1}.
\end{equation*}%
Hence
\begin{equation*}
\norm{ \frac{d\boldsymbol{B}_{\alpha }}{dt}} _{L^{2}\left( \left[
0,T\right] ;D\left( A\right) ^{\prime }\right) }\leq K,
\end{equation*}%
where $K$ is independent of $\alpha $.

From \eqref{eq:alphaMHD:Projected:velocity} we have%
\begin{equation*}
\norm{ A^{-1}\frac{d\boldsymbol{u}_{\alpha }}{dt}} \leq \norm{
A^{-1}\left( I+\alpha ^{2}A\right) ^{-1}\tilde{B}\left(
\boldsymbol{u},\boldsymbol{v}\right) } +\nu \norm{
\boldsymbol{u}_{\alpha }} +\norm{ A^{-1}\left( I+\alpha ^{2}A\right)
^{-1}B\left( \boldsymbol{B},\boldsymbol{B}\right) } ,
\end{equation*}%
and using \eqref{eq:Btilde_estimate_V_H_D(A)_short}%
\begin{align*}
\lnorm{ A^{-1}\left( I+\alpha ^{2}A\right) ^{-1}\tilde{B}\left(
\boldsymbol{u}_{\alpha },\boldsymbol{v}_{\alpha }\right) } & \leq
\lnorm{ A^{-1}\tilde{B}\left( \boldsymbol{u}_{\alpha },\boldsymbol{v}%
_{\alpha }\right) }  \\
& \leq c\lambda _{1}^{-1/4}\vnorm{ \boldsymbol{u}_{\alpha }}
\lnorm{ \boldsymbol{v}_{\alpha }}  \\
& \leq c\lambda _{1}^{-1/4}\vnorm{ \boldsymbol{u}_{\alpha }} \left(
\lnorm{ \boldsymbol{u}_{\alpha }} +\alpha ^{2}\lnorm{
A\boldsymbol{u}_{\alpha }} \right) ,
\end{align*}%
thus%
\begin{equation*}
\lnorm{ A^{-1}\left( I+\alpha ^{2}A\right) ^{-1}\tilde{B}\left(
\boldsymbol{u}_{\alpha },\boldsymbol{v}_{\alpha }\right) }
^{2}\leq 2c\lambda _{1}^{-1/2}{\tilde{k}_{1}}\left( \vnorm{ \boldsymbol{u}%
_{\alpha }} ^{2}+\alpha ^{2}\lnorm{ A\boldsymbol{u}_{\alpha }}
^{2}\right) ,
\end{equation*}%
and%
\begin{equation*}
{\int_{0}^{T}}\lnorm{ A^{-1}\left( I+\alpha ^{2}A\right) ^{-1}\tilde{B}%
\left( \boldsymbol{u}_{\alpha }\left( t\right)
,\boldsymbol{v}_{\alpha
}\left( t\right) \right) } ^{2}dt\leq c\lambda _{1}^{-1/2}{%
\tilde{k}_{1}^{2}}\nu ^{-1}.
\end{equation*}%
Also by \eqref{eq:BandBtilde_estimate_H_V_D(A)}%
\begin{equation*}
\norm{ B\left( \boldsymbol{B}_{\alpha },\boldsymbol{B}_{\alpha
}\right) } _{L^{2}\left( \left[ 0,T\right] ;D\left( A\right)
^{\prime }\right) }\leq c\lambda _{1}^{-1/2}{\tilde{k}_{1}^{2}\eta
}^{-1}.
\end{equation*}%
As a result we have%
\begin{equation*}
\norm{ \frac{d\boldsymbol{u}_{\alpha }}{dt}} _{L^{2}\left( \left[
0,T\right] ;D\left( A\right) ^{\prime }\right) }\leq K.
\end{equation*}%
Using Aubin's Compactness Lemma (see, for example, \cite[Lemma
8.4]{b_CF88})
we can extract subsequences of $\boldsymbol{u}_{\alpha _{j}}$ and $%
\boldsymbol{B}_{\alpha _{j}}$, which we relabel by
$\boldsymbol{u}_{\alpha
_{j}}$ and $\boldsymbol{B}_{\alpha _{j}}$ respectively, such that ${%
\boldsymbol{u}_{\alpha _{j}}}\rightarrow \boldsymbol{u}$ and ${\boldsymbol{B}%
_{\alpha _{j}}}\rightarrow \boldsymbol{B}$ strongly in $L^{2}\left(
\left[ 0,T\right] ;H\right) $ and strongly in $C\left( \left[
0,T\right] ;D\left( A\right) ^{\prime }\right) $, as $\alpha
_{j}\rightarrow 0^{+}$.

Observing that
\begin{equation*}
\norm{ \boldsymbol{v}_{\alpha _{j}}-\boldsymbol{u}_{\alpha _{j}}}
_{L^{2}\left( \left[ 0,T\right] ;V^{\prime }\right) }=\alpha
_{j}^{2}{\int_{0}^{T}}\vnorm{ \boldsymbol{u}_{\alpha _{j}}\left(
t\right) } dt\leq \alpha _{j}^{2}\nu ^{-1}{\tilde{k}_{1}},
\end{equation*}%
we obtain that $\boldsymbol{v}_{\alpha _{j}}\rightarrow \boldsymbol{u}$ in $%
L^{2}\left( \left[ 0,T\right] ;V^{\prime }\right) $, as $\alpha
_{j}\rightarrow 0^{+}$; and hence also that $\boldsymbol{u}\left( t\right) =%
\boldsymbol{v}\left( t\right) $ almost everywhere on $\left[
0,T\right] $.

Now, following the lines of the proof of Theorem
\ref{thm:alphaMHD:weakSol},
we can extract further subsequences (which we relabel by $\boldsymbol{u}%
_{\alpha _{j}},\boldsymbol{v}_{\alpha _{j}}$ and
$\boldsymbol{B}_{\alpha _{j}}$) and show that as $\alpha
_{j}\rightarrow 0^{+}$,
\begin{equation*}
\tilde{B}\left( \boldsymbol{u}_{\alpha _{j}},\boldsymbol{v}_{\alpha
_{j}}\right) \rightarrow \tilde{B}\left( \boldsymbol{v},\boldsymbol{v}%
\right) =B\left( \boldsymbol{v},\boldsymbol{v}\right)
\end{equation*}%
weakly in $L^{1}\left( \left[ 0,T\right] ;D\left( A\right) ^{\prime
}\right) $, and
\begin{align*}
B\left( \boldsymbol{B}_{\alpha _{j}},\boldsymbol{B}_{\alpha
_{j}}\right) \rightarrow B\left(
\boldsymbol{B},\boldsymbol{B}\right) ,\, %
B\left( \boldsymbol{u}_{\alpha _{j}},\boldsymbol{B}_{\alpha
_{j}}\right)
\rightarrow B\left( \boldsymbol{v},\boldsymbol{B}\right) ,\, %
B\left( \boldsymbol{B}_{\alpha _{j}},\boldsymbol{u}_{\alpha
_{j}}\right) \rightarrow B\left(
\boldsymbol{B},\boldsymbol{v}\right)
\end{align*}%
weakly in $L^{1}\left( \left[ 0,T\right] ;V^{\prime }\right) $.
Hence, we can pass to the limit (in the interpretation given by
\eqref{grp:alphaMHD:weakSol_integralFormulation}) in
\begin{align*}
& \left\langle \frac{d}{dt}\boldsymbol{\boldsymbol{v}}_{\alpha _{j}},%
\boldsymbol{w}\right\rangle _{D\left( A\right) ^{\prime
}}+\left\langle
\tilde{B}\left( \boldsymbol{u}_{\alpha _{j}},\boldsymbol{\boldsymbol{v}}%
_{\alpha _{j}}\right) ,\boldsymbol{w}\right\rangle _{D\left(
A\right)^{\prime }}+\nu \left( \boldsymbol{v}_{\alpha
_{j}},A\boldsymbol{w}\right) =\left\langle B\left(
\boldsymbol{B}_{\alpha _{j}},\boldsymbol{B}_{\alpha _{j}}\right) ,%
\boldsymbol{w}\right\rangle _{V^{\prime }}, \\
& \left\langle \frac{d}{dt}\boldsymbol{B}_{\alpha _{j}},\boldsymbol{\xi }%
\right\rangle _{V^{\prime }}+\left( B\left( \boldsymbol{u}_{\alpha _{j}},%
\boldsymbol{B}_{\alpha _{j}}\right) ,\boldsymbol{\xi }\right)
-\left(
B\left( \boldsymbol{B}_{\alpha _{j}},\boldsymbol{u}_{\alpha _{j}}\right) ,%
\boldsymbol{\xi }\right) +\eta \left( \left( \boldsymbol{B}_{\alpha _{j}},%
\boldsymbol{\xi }\right) \right) =0,
\end{align*}%
$\boldsymbol{w}\in D\left( A\right) ,\ \boldsymbol{\xi }\in V$  and
we
obtain that%
\begin{align*}
& \left\langle \frac{d}{dt}\boldsymbol{\boldsymbol{v}},\boldsymbol{w}%
\right\rangle _{D\left( A\right)  ^{\prime} }+\left\langle B\left( \boldsymbol{v}%
,\boldsymbol{\boldsymbol{v}}\right) ,\boldsymbol{w}\right\rangle
_{D(A) ^{\prime} }+\nu \left( \left(
\boldsymbol{v},\boldsymbol{w}\right)
\right) =\left\langle B\left( \boldsymbol{B},\boldsymbol{B}\right) ,%
\boldsymbol{w}\right\rangle _{V ^{\prime} }, \\
& \left\langle \frac{d}{dt}\boldsymbol{B},\boldsymbol{\xi
}\right\rangle
_{V^{\prime }}+\left( B\left( \boldsymbol{v},\boldsymbol{B}\right) ,%
\boldsymbol{\xi }\right) -\left( B\left( \boldsymbol{B},\boldsymbol{v}%
\right) ,\boldsymbol{\xi }\right) +\eta \left( \left( \boldsymbol{B},%
\boldsymbol{\xi }\right) \right) =0,
\end{align*}%
for every $\boldsymbol{w}\in D\left( A\right) ,\ \boldsymbol{\xi
}\in V$ and for almost every $t\in \left[ 0,T\right] $.

Now, since $\boldsymbol{v}\in L^{2}\left( \left[ 0,T\right]
;V\right) $, one can show that $B\left(
\boldsymbol{v},\boldsymbol{\boldsymbol{v}}\right) \in L^{1}\left(
\left[ 0,T\right] ;V^{\prime }\right) $ and then also that \mbox{$
\left({d}/{dt}\right)\boldsymbol{\boldsymbol{v}}\in L^{1}\left(
\left[ 0,T\right] ;V^{\prime }\right) $}, and since $w\in D\left(
A\right) $, which is dense in $V$, we
obtain the weak formulation of the MHD equations%
\begin{align*}
& \left\langle \frac{d}{dt}{\boldsymbol{v}},\boldsymbol{w}
\right\rangle _{V^{\prime }}+\left\langle B\left( \boldsymbol{v},\boldsymbol{%
\boldsymbol{v}}\right) ,\boldsymbol{w}\right\rangle _{V^{\prime
}}+\nu \left( \left( \boldsymbol{v},\boldsymbol{w}\right) \right)
=\left\langle B\left( \boldsymbol{B},\boldsymbol{B}\right)
,\boldsymbol{w}\right\rangle
_{V^{\prime }}, \\
& \left\langle \frac{d}{dt}\boldsymbol{B},\boldsymbol{\xi
}\right\rangle
_{V^{\prime }}+\left( B\left( \boldsymbol{v},\boldsymbol{B}\right) ,%
\boldsymbol{\xi }\right) -\left( B\left( \boldsymbol{B},\boldsymbol{v}%
\right) ,\boldsymbol{\xi }\right) +\eta \left( \left( \boldsymbol{B},%
\boldsymbol{\xi }\right) \right) =0,
\end{align*}%
for every $\boldsymbol{w},\ \boldsymbol{\xi }\in V$ and for almost every $%
t\in \left[ 0,T\right] $.

We notice, that 
every weak solution of \eqref{grp:alphaMHD:Projected} satisfies
the energy equality \eqref{eq:alphaMHD:weakSol_energyEquality} and
hence  the energy inequality \eqref{eq:MHD:energyIneq} follows by
passing to the $\liminf$ as $\alpha\to 0^{+}$, using the fact that
if $x_\alpha\to x$ weakly in a Hilbert space $X$, then $\|x\| \leq
\liminf\|x_\alpha\|$.
\end{proof}

\begin{theorem}
Let $T>0$, $\boldsymbol{u}^{in}\in D(A),\,\boldsymbol{B}^{in}\in V$
and denote by $\boldsymbol{u}_{\alpha },\,\boldsymbol{B}_{\alpha }$
and \mbox{$
\boldsymbol{v}_{\alpha }=\boldsymbol{u}_{\alpha }+\alpha ^{2}A\boldsymbol{u}%
_{\alpha }$} the strong solution of \eqref{grp:alphaMHD:Projected}
on $ \left[ 0,T\right] $. Then there exist $T_*=T_*(\Omega,  \nu,
\eta,  u^{in}, B^{in})$, $0<T_*\leq T$, subsequences
$\boldsymbol{u}_{\alpha _{j}},\,\boldsymbol{v}_{\alpha
_{j}},\,\boldsymbol{B}_{\alpha _{j}}$ and a pair of functions
$\boldsymbol{v},\boldsymbol{B}\in L^{\infty }\left( \left[
0,T_*\right] ;V\right) \cap L^{2}\left( \left[ 0,T_*\right]
;D(A)\right) $ such that, as \mbox{$\alpha _{j}\rightarrow 0^{+}$},

\begin{enumerate}
\item $\boldsymbol{u}_{\alpha _{j}}\rightarrow \boldsymbol{v}$ and $%
\boldsymbol{B}_{\alpha _{j}}\rightarrow \boldsymbol{B}$ weakly in $%
L^{2}\left( \left[ 0,T_*\right] ;D(A)\right) $ and strongly in
$L^{2}\left( \left[ 0,T_*\right] ;V\right) $,

\item $\boldsymbol{v}_{\alpha _{j}}\rightarrow \boldsymbol{v}$ weakly in $%
L^{2}\left( \left[ 0,T_*\right] ;V\right) $ and strongly in
$L^{2}\left( \left[ 0,T_*\right] ;H\right) $ and

\item $\boldsymbol{u}_{\alpha _{j}}\left( t\right) \rightarrow \boldsymbol{v}\left( t\right)
$ and $\boldsymbol{B}_{\alpha _{j}}\left( t\right) \rightarrow
\boldsymbol{B} \left( t\right)$ weakly in $V$ and uniformly on
$\left[ 0,T_*\right] $.
\end{enumerate}

Furthermore, the pair $\boldsymbol{v},\boldsymbol{B}$ is the unique
strong solution of the 3D MHD equations on $\left[0,T_*\right]$
with initial data $\boldsymbol{v}%
\left( 0\right) =\boldsymbol{u}^{in},\,\boldsymbol{B}\left( 0\right) =%
\boldsymbol{B}^{in}$. The strong solution of the 3D MHD equations
satisfies the energy equality
\begin{equation*}
\lnorm{ \boldsymbol{v}\left( t\right) } ^{2}+\lnorm{ \boldsymbol{B}
\left( t\right) } ^{2} +2 \int_{t_0}^{t}\left( \nu\vnorm{
\boldsymbol{v}(s)} ^{2}+\eta\vnorm{ \boldsymbol{B} (s)}
^{2}\right)ds %
= \lnorm{ \boldsymbol{v}\left( t_0\right) } ^{2}+\lnorm{
\boldsymbol{B}\left( t_0\right) } ^{2},\qquad  0\leq t_0\leq t \leq
T_*.
\end{equation*}
\end{theorem}
\begin{proof}
To prove the theorem we need to show that there exists $T_*$ such
that we have a uniform (independent of $\alpha $) bound on
\begin{equation}
{\vnorm{ \boldsymbol{u}_{\alpha }\left( t\right) } ^{2}+\alpha
^{2}}\lnorm{ {A\boldsymbol{u}_{\alpha }\left( t\right) }} {%
^{2}+\vnorm{ \boldsymbol{B}_{\alpha }(t)} ^{2}}
\label{eq:cnvgStrSol:bnd1}
\end{equation}%
and
\begin{equation}
\int_{0}^{T_*}\left( \nu (\lnorm{ A\boldsymbol{u}_{\alpha }(t)}
^{2}+\alpha ^{2}\lnorm{ A^{3/2}\boldsymbol{u}_{\alpha }(t)}
^{2})+\eta \lnorm{ A\boldsymbol{B}_{\alpha }(t)} ^{2}\right) dt
\label{eq:cnvgStrSol:bnd2}
\end{equation}%
in $\left[ 0,T_*\right] $. Then we can continue similarly to the
proof of the previous theorem, appropriately smoothing the data and
replacing $T$ by $ T_*$. Next we derive the formal estimates on
\eqref{eq:cnvgStrSol:bnd1} and \eqref{eq:cnvgStrSol:bnd2} that can
be proved rigorously using the Galerkin approximation scheme and
then passing to the limit using the proof of Theorem
\ref{thm:alphaMHD:weakSol}.

Let us recall \eqref{eq:alphaMHD:sum:inner_product_A}. By
\eqref{eq:BandBtilde_estimate_DA_V_H} and several applications of
Young's inequality we bound
\begin{align*}
\abs{ \left( \tilde{B}\left( \boldsymbol{u}_{\alpha},\boldsymbol{v}_{\alpha}\right) ,A%
\boldsymbol{u}_{\alpha}\right) } & \leq c\vnorm{ \boldsymbol{u}_{\alpha}%
} ^{1/2}\lnorm{ A\boldsymbol{u}_{\alpha}} ^{1/2}(\vnorm{
\boldsymbol{u}_{\alpha}} +\alpha ^{2}\lnorm{ A^{3/2}\boldsymbol{u}_{\alpha}%
} )\lnorm{ A\boldsymbol{u}_{\alpha}}  \\
& \leq c\nu ^{-3}\vnorm{ \boldsymbol{u}_{\alpha}} ^{6}+\nu
^{-3}\alpha
^{6}\lnorm{ A\boldsymbol{u}_{\alpha}} ^{6}+\frac{\nu }{4}\lnorm{ A%
\boldsymbol{u}_{\alpha}} ^{2}+\frac{\nu }{2}\alpha ^{2}\lnorm{ A^{3/2}%
\boldsymbol{u}_{\alpha}} ^{2}.
\end{align*}%
By \eqref{eq:BandBtilde_estimate_DA_V_H}%
\begin{align*}
\abs{ \left( B\left( \boldsymbol{B}_{\alpha},\boldsymbol{B}_{\alpha}\right) ,A%
\boldsymbol{u}_{\alpha}\right) } & \leq c\vnorm{ \boldsymbol{B}_{\alpha}%
} ^{1/2}\lnorm{ A\boldsymbol{B}_{\alpha}} ^{1/2}\vnorm{
\boldsymbol{B}_{\alpha}} \lnorm{ A\boldsymbol{u}_{\alpha}}  \\
& \leq c\nu ^{-2}\eta ^{-1}\vnorm{ \boldsymbol{B}_{\alpha}} ^{6}+\frac{%
\eta }{4}\lnorm{ A\boldsymbol{B}_{\alpha}} ^{2}+\frac{\nu
}{4}\lnorm{ A\boldsymbol{u}_{\alpha}} ^{2}.
\end{align*}%
By \eqref{eq:BandBtilde_estimate_DA_V_H} we also have
\begin{align*}
\abs{ \left( B\left( \boldsymbol{B}_{\alpha},\boldsymbol{u}_{\alpha}\right) ,A%
\boldsymbol{B}_{\alpha}\right) } & \leq c\vnorm{ \boldsymbol{B}_{\alpha}%
} ^{1/2}\vnorm{ \boldsymbol{u}_{\alpha}} \lnorm{ A%
\boldsymbol{B}_{\alpha}} ^{3/2} \\
& \leq c\eta ^{-3}\vnorm{ \boldsymbol{B}_{\alpha}} ^{6}+\eta
^{-3}\vnorm{ \boldsymbol{u}_{\alpha}} ^{6}+\frac{\eta }{8}\lnorm{ A%
\boldsymbol{B}_{\alpha}} ^{2}
\end{align*}%
and by \eqref{eq:B_estimate_V_DA_H}
\begin{equation*}
\abs{ \left( B\left( \boldsymbol{u}_{\alpha},\boldsymbol{B}_{\alpha}\right) ,A%
\boldsymbol{B}_{\alpha}\right) } \leq c\vnorm{ \boldsymbol{B}_{\alpha}%
} ^{1/2}\vnorm{ \boldsymbol{u}_{\alpha}} \lnorm{ A%
\boldsymbol{B}_{\alpha}} ^{3/2}.
\end{equation*}%
Hence by \eqref{eq:alphaMHD:sum:inner_product_A} and the above
estimates we have
\begin{equation}
\frac{d}{dt}\left( \vnorm{ \boldsymbol{u}_{\alpha}} ^{2}+\alpha
^{2}\lnorm{ A\boldsymbol{u}_{\alpha}} ^{2}+\vnorm{ \boldsymbol{B}_{\alpha}%
} ^{2}\right) +\nu \left( |{A\boldsymbol{u}_{\alpha}|}^{2}+\alpha ^{2}|{%
A^{3/2}\boldsymbol{u}_{\alpha}|}^{2}\right) +\eta
|{A\boldsymbol{B}_{\alpha}|}^{2}\leq c\mu
^{-3}\left( \vnorm{ \boldsymbol{u}_{\alpha}} ^{6}+\alpha ^{6}\lnorm{ A%
\boldsymbol{u}_{\alpha}} ^{6}+\vnorm{ \boldsymbol{B}_{\alpha}}
^{6}\right) \label{eq:alphaMHD:ConvStrSols:H2_inequality}
\end{equation}%
Denote
\begin{equation*}
\boldsymbol{y}=\vnorm{ \boldsymbol{u}_{\alpha}} ^{2}+\alpha
^{2}\lnorm{ A\boldsymbol{u}_{\alpha}} ^{2}+\vnorm{
\boldsymbol{B}_{\alpha} } ^{2}.
\end{equation*}%
Now, if $y(0)=0$, that is
$\boldsymbol{u}^{in}=\boldsymbol{B}^{in}=0$, then the solution is
steady $\boldsymbol{u}_{\alpha}(t)\equiv 0$,
$\boldsymbol{v}_{\alpha}(t)\equiv 0$,
$\boldsymbol{B}_{\alpha}(t)\equiv 0$ and $\boldsymbol{v}(t)\equiv
0$, $\boldsymbol{B}(t)\equiv 0$ 
 exists for all
$t\geq 0$. Otherwise, from
\eqref{eq:alphaMHD:ConvStrSols:H2_inequality} we have
\begin{equation*}
\frac{d}{dt}\boldsymbol{y}\leq c\mu ^{-3}\boldsymbol{y}^{3}
\end{equation*}
and thus%
\begin{equation*}
\boldsymbol{y}\left( t\right) \leq 2\boldsymbol{y}\left( 0\right)
\end{equation*}%
for $0\leq t\leq \frac{3}{8}c^{3}\mu ^{3}\boldsymbol{y}\left(
0\right) ^{-2}$. We conclude that
\begin{equation*}
{\vnorm{ \boldsymbol{u}_{\alpha }\left( t\right) } ^{2}+\alpha
^{2}}\lnorm{ {A\boldsymbol{u}_{\alpha }\left( t\right) }} {%
^{2}+\vnorm{ \boldsymbol{B}_{\alpha }(t)} ^{2}}\leq 2\left(
\vnorm{ \boldsymbol{u}^{in}} ^{2}+\alpha ^{2}\lnorm{ A%
\boldsymbol{u}^{in}} ^{2}+\vnorm{ \boldsymbol{B}%
^{in}} ^{2}\right)
\end{equation*}%
for $0\leq t\leq T_*:=\min \left( T,\frac{3}{8}c^{3}\mu ^{3}\boldsymbol{y}%
\left( 0\right) ^{-2}\right) $. Also, by integrating %
\eqref{eq:alphaMHD:ConvStrSols:H2_inequality} over $(0,T_*)$, we obtain%
\begin{multline*}
\int_{0}^{T_*}\left( \nu (\lnorm{ A\boldsymbol{u}_{\alpha }(t)}
^{2}+\alpha ^{2}\lnorm{ A^{3/2}\boldsymbol{u}_{\alpha }(t)}
^{2})+\eta \lnorm{ A\boldsymbol{B}_{\alpha
}(t)} ^{2}\right) dt \\
\leq \vnorm{ \boldsymbol{u}^{in}} ^{2}+\alpha ^{2}\lnorm{ A%
\boldsymbol{u}^{in}} ^{2}+\vnorm{ \boldsymbol{B}%
^{in}} ^{2}+c\mu ^{-3}T_*\left( \vnorm{ \boldsymbol{u}%
^{in}} ^{2}+\alpha ^{2}\lnorm{ A\boldsymbol{u}^{in}} ^{2}+\vnorm{
\boldsymbol{B}^{in}} ^{2}\right) ^{3}.
\end{multline*}%
Assuming that $0<\alpha \leq L$, the bounds are independent of
$\alpha $.
\end{proof}

\section{Discussion}
We proved the well-posedness of the three-dimensional MHD-$\alpha$
model \eqref{grp:alphaMHD_intro} in the periodic boundary
conditions. This model modifies the nonlinearity of the MHD
equations \eqref{grp:MHD} without enhancing dissipation. We showed
that the model has a unique global weak (or strong, for smooth
initial data) solution. Also, there is a subsequence of weak
solutions of the MHD-$ \alpha $ equations that converge, as
\mbox{$\alpha \rightarrow $ $0^{+}$}, (in the appropriate sense) to
a Leray-Hopf weak solution (which satisfies the energy inequality
\eqref{eq:MHD:energyIneq}) of the MHD equations \eqref{grp:MHD} on
any time interval $[0,T]$. Also, if the initial data is smooth, a
subsequence of solutions converges for a short interval of time, to
the unique strong solution of the MHD equations on this interval.
These properties are essential for the $\alpha$ models to be
regarded as regularizing numerical schemes. In a follow up paper, we
intend to do the error estimates in which we will investigate the
error in terms of $m$ and $\alpha$. Namely, the distance between the
solution of the Galerkin MHD-$\alpha$ model to that of the exact
strong solution of the MHD equations, for smooth initial data.

There are many different $\alpha$ models. For example, the global
well-posedness can be shown for the 3D Modified-Leray-$ \alpha $-MHD
model \eqref{grp:ML_alpha_MHD_intro}. However, at the moment we are
unable to find a conserved quantity in the ideal version of
\eqref{grp:ML_alpha_MHD_intro}, which can be identified with a cross
helicity, contrary to the MHD-$\alpha$ model
\eqref{grp:alphaMHD_intro}, where there exist the ideal invariants
that could be identified with the three invariants of the original
3D MHD equations.

\section*{Acknowledgements}
E.S.T. is thankful for the kind hospitality of the  \'{E}cole
Normale Sup\'{e}rieure - Paris where this work was completed. This
work was supported in part by the NSF, grant no.~DMS-0204794,  the
ISF Grant no.~120/06, and the BSF, grant no.~200423.

\bibliographystyle{amsplain}
\bibliography{MHDBib}

\end{document}